\documentclass[12pt]{article}
\usepackage{color}
\definecolor{darkblue}{rgb}{0.00,0.25,0.50}
\usepackage[colorlinks,filecolor=blue,citecolor=darkblue]{hyperref}

\usepackage{amsfonts,amssymb,amsmath,amsthm}
\usepackage{url}
\usepackage{enumerate}
\usepackage[ukrainian, russian, english]{babel}
\usepackage[cp1251]{inputenc}

\usepackage{cite}
\begin{document}
\selectlanguage{ukrainian} 
\thispagestyle{empty}

\title{}

\begin{center}
\textbf{\Large Точні оцінки колмогоровських поперечників класів аналітичних функцій}
\end{center}
\vskip0.5cm
\begin{center}
А.\,С.~Сердюк, В.\,В.~Боденчук\\ \emph{\small Інститут математики НАН України, Київ}
\end{center}
\vskip0.5cm

\begin{abstract}
We prove that kernels of analytic functions of kind ${H_{h,\beta}(t)=\sum\limits_{k=1}^{\infty}\frac{1}{\ch kh}\cos\Big(kt-\frac{\beta\pi}{2}\Big)}$, $h>0$, ${\beta\in\mathbb{R}}$, satisfies Kushpel's condition $C_{y,2n}$ beginning with some number $n_h$ which is explicitly expressed by parameter $h$ of smoothness of the kernel.
As a consequence, for all $n\geqslant n_h$ we obtain lower bounds for Kolmogorov widths $d_{2n}$ of functional classes that are representable as convolutions of kernel $H_{h,\beta}$ with functions $\varphi\perp1$, which belong to the unit ball in the space $L_{\infty}$, in the space $C$.
The obtained estimates coincide with the best uniform approximations by trigonometric polynomials for these classes. As a result, we obtain exact values for widths of mentioned classes of convolutions.
Also for all $n\geqslant n_h$ we obtain exact values for Kolmogorov widths $d_{2n-1}$ of classes of convolutions of functions $\varphi\perp1$, which belong to the unit ball in the space $L_1$, with kernel $H_{h,\beta}$ in the space $L_1$.

\vskip0.6cm

У роботі встановлено, що ядра аналітичних функцій вигляду ${H_{h,\beta}(t)=\sum\limits_{k=1}^{\infty}\frac{1}{\ch kh}\cos\Big(kt-\frac{\beta\pi}{2}\Big)}$, $h>0$, ${\beta\in\mathbb{R}}$, задовольняють введену Кушпелем умову $C_{y,2n}$, починаючи з деякого номера $n_h$, який в явному вигляді виражається через параметр $h$ гладкості ядра.
В результаті отримано оцінки знизу колмогоровських поперечників $d_{2n}$ для всіх $n\geqslant n_h$ функціональних класів, що зображуються згортками ядра $H_{h,\beta}$ із функціями $\varphi\perp1$, що належать одиничній кулі простору $L_{\infty}$, в просторі $C$. Отримані оцінки співпали з найкращими рівномірними наближеннями зазначених класів тригонометричними поліномами. Як наслідок, знайдено точні значення поперечників вказаних класів згорток.
Знайдено також точні значення поперечників $d_{2n-1}$ в просторі $L_1$ для усіх $n\geqslant n_h$ класів згорток функцій $\varphi\perp1$, що належать одиничній кулі простору $L_1$, з ядром $H_{h,\beta}$.
\end{abstract}

\vskip1cm


\textbf{1. Постановка задачі.} Через $L=L_1$ позначимо простір $2\pi$-періодичних сумовних функцій $f$ з нормою
$\|f\|_1=\int\limits_{-\pi}^{\pi}|f(t)|dt$, через
$L_\infty$ --- простір $2\pi$-періодичних вимірних і суттєво обмежених функцій з нормою
${\|f\|_\infty=\mathop{\rm ess\,sup}\limits_{t\in\mathbb{R}\ }|f(t)|,}$
а через $C$ --- простір $2\pi$-періодичних неперервних функцій $f$, у якому норма задається рівністю
$\|f\|_C=\max\limits_{t\in\mathbb{R}}|f(t)|$.

Нехай $\Psi_\beta(t)$ --- фіксоване сумовне ядро вигляду
\begin{equation}\label{Psi_beta}
\Psi_\beta(t)=\sum\limits_{k=1}^{\infty}\psi(k)\cos\left(kt-\frac{\beta\pi}{2}\right),\,\psi(k)>0,\, \sum\limits_{k=1}^{\infty}\psi(k)<\infty,\, \beta\in\mathbb{R}.
\end{equation}
Через $C_{\beta,p}^{\psi}$, $p=1, \infty$, позначимо клас $2\pi$-періодичних функцій $f$, що зображуються у вигляді згортки з ядром $\Psi_\beta$
\begin{equation}\label{f}
f(x)
	=
		A+\left(\Psi_\beta\ast\varphi\right)(x)
	=
		A+\frac{1}{\pi}\int\limits_{-\pi}^{\pi}\Psi_\beta(x-t)\varphi(t)dt,
		\; A\in\mathbb{R},
\end{equation}
де $\|\varphi\|_p\leqslant1,\, \varphi\perp1.$
Функцію $\varphi$ в рівності \eqref{f} називають $(\psi,\beta)$-похідною функції $f$ і позначають через $f_\beta^\psi$.
Поняття $(\psi,\beta)$-похідної введене О.І. Степанцем (див., наприклад, \cite[\S~7--8]{Stepanets_2002_1}).

У роботі розглядаються ядра $\Psi_\beta(t)$ вигляду \eqref{Psi_beta} при $\psi(k)=\frac{1}{\ch kh}$, $h>0$, тобто функції
\begin{equation}\label{Psi_kernel}
H_{h,\beta}(t)
	=
		\sum\limits_{k=1}^{\infty}
			\frac{1}{\ch kh}
			\cos
				\Big(
					kt-\frac{\beta\pi}{2}
				\Big),\, 
	h>0,\, 
	\beta\in\mathbb{R}.
\end{equation}
При зазначених $\psi$ класи $C_{\beta,p}^{\psi}$ будемо позначати через $C_{\beta,p}^{h}$.

Як показано в \cite[с.~141]{Stepanets_2002_1}, функції з класів $C_{\beta,p}^{h}$, $h>0$, складаються із функцій $f\in C$, що допускають регулярне продовження $f(z)=f(x+iy)=u(x,y)+iv(x,y)$ у смугу
\begin{equation}\label{strip}
\{z=x+iy:\, -h<y<h\}
\end{equation}
комплексної площини. Зокрема (див. \cite[с.~269]{Akhiezer_1965}), при $\beta=2l$, $l\in\mathbb{Z}$, і $p=\infty$ класи $C_{\beta,p}^{h}$ збігаються з відомими класами $A_{\infty}^{h}$ фукнцій $f\in C$, які допускають аналітичне продовження у смугу \eqref{strip} і таких, що $\|\mathop{\textrm{Re}}f(\cdot+iy)\|_{\infty}\leqslant1$, $|y|<h$.

Нехай $d_m(\mathfrak{N},X)$ --- поперечник за Колмогоровим порядку $m$ центрально-симетричної множини $\mathfrak{N}\subset X$ в банаховому просторі $X$, тобто величина вигляду
\begin{equation}\label{d_m}
d_m(\mathfrak{N},X)
	=
		\inf_{F_m\subset X}
		\sup_{f\in \mathfrak{N}}
		\inf_{y\in F_m}\|f-y\|_X,
\end{equation}
де зовнішній $\inf$ розглядається по всіх $m$-вимірних лінійних підпросторах $F_m$ із $X$.

Розв’язується задача знаходження точних значень поперечників $d_{2n}(C_{\beta,\infty}^{h},C)$, $d_{2n-1}(C_{\beta,\infty}^{h},C)$ та $d_{2n-1}(C_{\beta,1}^{h},L)$ для довільних $h>0$, $\beta\in\mathbb{R}$ та усіх натуральних $n$, більших деякого номера, залежного лише від параметра $h$.

Позначивши через $\mathcal{T}_{2n-1}$ підпростір тригонометричних поліномів $t_{n-1}$ порядку $n-1$, розглянемо величини найкращих наближень
\begin{equation}\label{E_n}
E_n(\mathfrak{N})_X
	=
		\sup_{f\in \mathfrak{N}}
		\inf_{t_{n-1}\in \mathcal{T}_{2n-1}}\|f-t_{n-1}\|_X,
		\; \mathfrak{N}\subset X\subset L.
\end{equation}

Із \eqref{d_m} і \eqref{E_n} випливає, що при всіх $n\in\mathbb{N}$
\begin{equation}\label{d_m_E_n}
d_{2n-1}(\mathfrak{N},X)
	\leqslant 
		E_n(\mathfrak{N})_X,
		\; \mathfrak{N}\subset X\subset L.
\end{equation}

З робіт Н.І. Ахієзера \cite{Akhiezer_1938} та С.М. Нікольського \cite{Nikolskiy_1946} випливає, що при $\beta=2l$, $l\in\mathbb{Z}$, і довільних $h>0$ та $n\in\mathbb{N}$ виконуються рівності
\begin{equation*}
E_n(C_{\beta,\infty}^{h})_C
	=
		E_n(C_{\beta,1}^{h})_{L}
	=
		\|H_{h,\beta}\ast\varphi_n\|_C
	=
		\frac{4}{\pi}
		\sum\limits_{\nu=0}^{\infty}
			\frac
				{(-1)^{\nu}}
				{(2\nu+1)\ch((2\nu+1)nh)},
\end{equation*}
де
\begin{equation}\label{varp}
\varphi_n(t):=\textnormal{sign}\sin nt.
\end{equation}

В даній роботі ми покажемо, що при довільних $\beta\in\mathbb{R}$ точні значення величин $E_n(C_{\beta,\infty}^{h})_C$ та $E_n(C_{\beta,1}^{h})_{L}$ для усіх номерів $n$, починаючи з деякого номера $n_h^*$, можна одержати на основі результатів роботи \cite{Serdyuk_2002}. При цьому виявляється, що для усіх зазначених $n$ рівності
\begin{equation*}
E_n(C_{\beta,\infty}^{h})_C
	=
		E_n(C_{\beta,1}^{h})_{L}
	=
		\|H_{h,\beta}\ast\varphi_n\|_C
\end{equation*}
залишаються вірними при довільних $\beta\in\mathbb{R}$.

Тому для розв'язання задачі про знаходження точних значень колмогоровських поперечників залишається встановити оцінки знизу
\begin{equation}\label{dno1}
d_{2n}(C_{\beta,\infty}^{h}, C)
	\geqslant
		\|H_{h,\beta}\ast\varphi_n\|_C,
\end{equation}
\begin{equation}\label{dno2}
d_{2n-1}(C_{\beta,1}^{h}, L)
	\geqslant
		\|H_{h,\beta}\ast\varphi_n\|_C.
\end{equation}
для усіх номерів $n$, починаючи з деякого номера $n_h$.

При $\beta=2l$, $l\in\mathbb{Z}$, В.М.~Тихомиров у \cite{Tihomirov_1960} та \cite{Tihomirov_1976} одержав рівності
\begin{equation}\label{TihRivn}
d_{2n-1}(C_{\beta,\infty}^{h}, C)
	=
		d_{2n}(C_{\beta,\infty}^{h}, C)
	=
		\|H_{h,\beta}\ast\varphi_n\|_C, n\in\mathbb{N}.
\end{equation}
Однак доведення рівностей \eqref{TihRivn} у вказаних роботах не були повними. Коректне доведення зрештою було отримане Форстом \cite{Forst_1977}, який фактично показав, що ядро $H_{h,\beta}(t)$ при $\beta=2l$, $l\in\mathbb{Z}$, є $\text{CVD}$-ядром (ядром, що не збільшує осциляції). Кажуть, що ядро $K\in L$ є $\text{CVD}$-ядром і записують $K\in \text{CVD}$, якщо для довільної функції $\varphi\in C$ виконується нерівність $\nu(K\ast \varphi)\leqslant\nu(\varphi)$, де $\nu(g)$ --- число змін знаку на $[0,2\pi)$ функції ${g\in C}$. Згодом А.~Пінкус розробив методи, які дозволяють отримувати точні оцінки поперечників для класів згорток, що породжуються довільними $\text{CVD}$-ядрами (див. \cite{Pinkus_1985, Pinkus_1979}).
Якщо ж $\beta\not=2l$, $l\in\mathbb{Z}$, то, як показано в \cite[с.~111]{Kushpel_1988_dys}, ядра $H_{h,\beta}(t)$ можуть збільшувати осциляції, тому точні оцінки знизу поперечників $d_{2n}(C_{\beta,\infty}^{h},C)$ та $d_{2n-1}(C_{\beta,1}^{h},L)$ неможливо отримати, користуючись методами, які розвинуто А. Пінкусом \cite{Pinkus_1985}.

Зауважимо також, що для усіх $h>0$ таких, що
\begin{equation}\label{relatio111}
\frac{\ch kh}{\ch (k+1)h}
	\leqslant
		\frac{\ch h}{\ch 2h}
	\leqslant
		\rho(\beta),\; k\in\mathbb{N},
\end{equation}
де $\rho(\beta)=0{,}2$, якщо $\beta\in\mathbb{Z}$ і $\rho(\beta)=0{,}193864$, якщо ${\beta\in\mathbb{R\setminus Z}}$
нерівності \eqref{dno1} та \eqref{dno2} при довільних $n\in\mathbb{N}$ випливають з роботи \cite[с.~1118--1119]{Serdyuk_1995}. Обчислення показують, що умова \eqref{relatio111}, а разом з нею і оцінки \eqref{dno1} та \eqref{dno2}, має місце при усіх $h\geqslant 1{,}644651$, якщо $\beta\in\mathbb{Z}$ і $h\geqslant 1{,}67423$, якщо ${\beta\in\mathbb{R\setminus Z}}$.

\textbf{2. Формулювання основних результатів.} 
Для кожного фіксованого $h>0$ покладемо
\begin{equation*}
n_h^*=
\begin{cases}
1,& \text{якщо $h\geqslant\ln\frac{10}{3}$,}
\\
n_h^{**},& \text{якщо $0<h<\ln\frac{10}{3}$,}
\end{cases}
\end{equation*}
де $n_h^{**}$ --- найменше натуральне число, для якого виконується нерівність
\begin{equation}\label{umova_n0_h}
(1-e^{-h})^2
	\geqslant
		\frac{5+3e^{-2h}}{1-e^{-2h}}
        \frac
            {\left(\frac{1+e^{-2h}}{2}\right)^{2n}}
            {\sqrt{1-\left(\frac{1+e^{-2h}}{2}\right)^{2n}}}		
		+
		(2+e^{-2nh})e^{-2nh}.
\end{equation}
Має місце твердження.

\textbf{Теорема 1.} \textit{Нехай $h>0$ і $\beta\in\mathbb{R}$. 
Тоді для всіх номерів $n$ таких, що $n\geqslant n_h^*$, виконуються рівності
\begin{equation*}
E_n(C_{\beta,\infty}^{h})_C
	=
		E_n(C_{\beta,1}^{h})_{L}
	=
		\|H_{h,\beta}\ast\varphi_n\|_C
	=
\end{equation*}
\begin{equation}\label{E_n_rivnosti}
	=
	\frac{4}{\pi}
	\left|
		\sum\limits_{\nu=0}^{\infty}
			\frac{1}{(2\nu+1)\ch((2\nu+1)nh)}
			\sin
			\left(
				(2\nu+1)\theta_n\pi-\frac{\beta\pi}{2}
			\right)
	\right|,
\end{equation}
в яких $\varphi_n(t)$ --- функція вигляду \eqref{varp}, а $\theta_n=\theta_n(h,\beta)$ --- єдиний на $[0,1)$ корінь рівняння}
\begin{equation}\label{theta}
\sum\limits_{\nu=0}^{\infty}
	\frac{1}{\ch((2\nu+1)nh)}
	\cos
	\left(
		(2\nu+1)\theta_n\pi-\frac{\beta\pi}{2}
	\right)
	=0.
\end{equation}

Далі для кожного фіксованого $h>0$ через $n_{h}$ будемо позначати найменший з номерів $n\geqslant9$, для якого виконується нерівність
\begin{equation*}
\frac{37}{5(1-e^{-h})}e^{-h\sqrt{n}}
+\frac{e^{-h}}{(1-e^{-h})^2} 
	\min\left\{\frac{160}{27(n-\sqrt{n})}, \frac{8}{3n-7\sqrt{n}}\right\}
		\leqslant
\end{equation*}
\begin{equation}\label{umova_n_0}
		\leqslant
			\left(\frac{1}{2}+\frac{1}{(1-e^{-h})\ch h}\right)
			\left(\frac{1-e^{-h}}{1+e^{-h}}\right)^{\frac {4}{1-e^{-2h}}}
\end{equation}

У прийнятих позначеннях мають місце наступні твердження.

\textbf{Теорема 2.} \textit{Нехай $h>0$, $\beta\in \mathbb{R}$.
Тоді для всіх номерів $n$ таких, що $n\geqslant n_{h}$, виконуються нерівності \eqref{dno1} і \eqref{dno2}}.

\textbf{Теорема 3.} \textit{Нехай $h>0$, $\beta\in \mathbb{R}$.
Тоді для всіх номерів $n$ таких, що $n\geqslant n_{h}$, виконуються рівності }
\begin{equation*}
d_{2n}(C_{\beta,\infty}^{h},C)
	=
		d_{2n-1}(C_{\beta,\infty}^{h},C)
	=
		d_{2n-1}(C_{\beta,1}^{h},L)
	=
\end{equation*}
\begin{equation*}
	=
		E_n(C_{\beta,\infty}^{h})_C
	=
		E_n(C_{\beta,1}^{h})_{L}
	=
		\|H_{h,\beta}\ast\varphi_n\|_C
	=
\end{equation*}
\begin{equation}\label{dn}
	=
	\frac{4}{\pi}
	\left|
		\sum\limits_{\nu=0}^{\infty}
			\frac
				{1}
				{(2\nu+1)\ch((2\nu+1)nh)}
			\sin
			\left(
				(2\nu+1)\theta_n\pi-\frac{\beta\pi}{2}
			\right)
	\right|,
\end{equation}
\textit{де $\theta_n=\theta_n(h,\beta)$ --- єдиний на $[0,1)$ корінь рівняння \eqref{theta}}.

При $\beta=2l-1$, $l\in\mathbb{Z}$, із теореми 3 одержуємо наступне твердження.

\textbf{Наслідок 1.}
\textit{Нехай $h>0$, $\beta=2l-1$, $l\in\mathbb{Z}$. Тоді для всіх номерів таких, що $n\geqslant n_{h}$, виконуються рівності}
\begin{equation*}
d_{2n}(C_{\beta,\infty}^{h},C)
	=
		d_{2n-1}(C_{\beta,\infty}^{h},C)
	=
		d_{2n-1}(C_{\beta,1}^{h},L)
	=
		E_n(C_{\beta,\infty}^{h})_C
	=
\end{equation*}
\begin{equation*}
	=
		E_n(C_{\beta,1}^{h})_{L}
	=
		\|H_{h,\beta}\ast\varphi_n\|_C
	=
		\frac{4}{\pi}
		\left|
			\sum\limits_{\nu=0}^{\infty}
				\frac
					{1}
					{(2\nu+1)\ch((2\nu+1)nh)}
		\right|.
\end{equation*}

Із теореми~3 легко одержати асимптотичні при ${n\to\infty}$ оцінки поперечників $d_{2n}(C_{\beta,\infty}^{h},C)$, $d_{2n-1}(C_{\beta,\infty}^{h},C)$ та $d_{2n-1}(C_{\beta,1}^{h},L)$.

\textbf{Теорема 4.}
\textit{Нехай $h>0$ та $\beta\in \mathbb{R}$.
Тоді при $n\geqslant n_{h}$ }
\begin{equation*}
d_{2n}(C_{\beta,\infty}^{h},C)
	=
		d_{2n-1}(C_{\beta,\infty}^{h},C)
	=
		d_{2n-1}(C_{\beta,1}^{h},L)
	=
		E_n(C_{\beta,\infty}^{h})_C=
\end{equation*}
\begin{equation}\label{Th_3}
	=
		E_n(C_{\beta,1}^{h})_L
	=
		\frac{1}{\ch nh}
		\left(
			\frac{4}{\pi}
			+\gamma_n\frac{e^{-2nh}}{1-e^{-2nh}}
		\right),
\end{equation}
\textit{де $|\gamma_n|\leqslant\frac{28}{3\pi}$.}

\textbf{3. Доведення теореми 1.} 
В роботі \cite{Serdyuk_2002} було встановлено, що якщо послідовність коефіцієнтів $\psi(k)$ ядра $\Psi_{\beta}$ вигляду \eqref{Psi_beta}, яке породжує класи $C_{\beta,p}^{\psi}$, $p=1,\infty$, задовольняє умову Даламбера
\begin{equation}\label{Dqumova}
\lim_{k\to\infty}\frac{\psi(k+1)}{\psi(k)}=q,\; q\in(0,1),\; \psi(k)>0,
\end{equation}
(при цьому записуватимемо: $\psi\in\mathcal{D}_q$), то знайдеться номер $n_0$ такий, що для будь-якого натурального $n\geqslant n_0$ мають місце рівності
\begin{equation*}
E_n(C_{\beta,\infty}^{\psi})_C
	=
		E_n(C_{\beta,1}^{\psi})_{L}
	=
		\|\Psi_{\beta}\ast\varphi_n\|_C
	=
\end{equation*}
\begin{equation}\label{E_n_rivnosti_psi}
	=
	\frac{4}{\pi}
	\left|
		\sum\limits_{\nu=0}^{\infty}
			\frac{\psi((2\nu+1)n)}{2\nu+1}
			\sin
			\left(
				(2\nu+1)\theta_n\pi-\frac{\beta\pi}{2}
			\right)
	\right|,
\end{equation}
в яких $\theta_n=\theta_n(\psi,\beta)$ --- єдиний на $[0,1)$ корінь рівняння
\begin{equation*}
\sum\limits_{\nu=0}^{\infty}
	\psi((2\nu+1)n)
	\cos
	\left(
		(2\nu+1)\theta_n\pi-\frac{\beta\pi}{2}
	\right)
	=0.
\end{equation*}
При цьому (див. \cite[с.~188--190]{Serdyuk_2002}) номер $n_0$ означається конструктивно як найменше натуральне число, для котрого виконуються нерівності
\begin{equation}\label{umova_n0}
(1-q)^2
	\geqslant
		\frac{5+3q^2}{1-q^2}
		\frac
			{\left(
				\frac{1+q^2}{2}
			\right)^{2n}}
			{\sqrt{1
				-\left(
					\frac{1+q^2}{2}
				\right)^{2n}}}
		+
		\varepsilon_n(2+\varepsilon_n),
		\; n=n_0, n_0+1, \dots,
\end{equation}
де
\begin{equation}\label{epsilonn}
\varepsilon_n
	=
		\varepsilon_n(\psi)
	:=
		\sup_{k\geqslant n}
			|\frac{\psi(k+1)}{\psi(k)}-q|.
\end{equation}

Крім того, згідно з теоремою~2 роботи \cite{Serdyuk_2002}, рівності \eqref{E_n_rivnosti_psi} справджуються для усіх $\beta\in\mathbb{R}$ і $n\in\mathbb{N}$ за умови, що $\frac{\psi(k+1)}{\psi(k)}<\rho_*,\; k=1,2,\dots,$ де ${\rho_*=0{,}3253678\dots}$ --- корінь рівняння
\begin{equation*}
2\rho
+
\frac
	{(1+3\rho)\rho^2}
	{(1-\rho)\sqrt{1-2\rho^2}}
		=1
\end{equation*}
на інтервалі (0,1).

Оскільки коефіцієнти $\psi(k)=\frac{1}{\ch kh}$ ядра $H_{h,\beta}(t)$ задовольняють умову $\mathcal{D}_q$ при $q=e^{-h}$ і відношення $\frac{\psi(k+1)}{\psi(k)}=\frac{\ch kh}{\ch (k+1)h}=q\frac{1+q^{2k}}{1+q^{2k+2}}$ утворює спадну послідовність, то при  $q\in(0,\frac{3}{10}]$
\begin{equation*}
\frac{\psi(k+1)}{\psi(k)}\leqslant q\frac{1+q^{2}}{1+q^{4}}\leqslant 0{,}3253678, k=1,2,\dots.
\end{equation*}
Отже, якщо $h>\ln\frac{10}{3}$, то рівності \eqref{E_n_rivnosti} виконуються для усіх $n\in\mathbb{N}$. При $\psi(k)=\frac{1}{\ch kh}$ для величин $\varepsilon_n=\varepsilon_n(\psi)$ вигляду \eqref{epsilonn} мають місце співвідношення
\begin{equation*}
\varepsilon_n
	=
		q^{2k+1}\frac{1-q^2}{1+q^{2k+2}}
	<
		q^{2k+1}, q=e^{-h}.
\end{equation*}
Тому виконання нерівності \eqref{umova_n0_h} гарантує виконання умови \eqref{umova_n0}, а отже, всилу \cite{Serdyuk_2002}, і рівностей \eqref{E_n_rivnosti} для усіх номерів $n$ таких, що $n\geqslant n_h^{**}$. Теорему доведено.

\textbf{4. Означення і допоміжні твердження.}
Доведення нерівностей \eqref{dno1} і \eqref{dno2} будемо здійснювати використовуючи запропонований О.К. Кушпелем \cite{Kushpel_1988} метод знаходження оцінок знизу поперечників класів згорток із твірними ядрами $\Psi_\beta$, що задовольняють так звану умову $C_{y,2n}$.
Наведемо означення і відомі твердження, які будуть використовуватись у подальшому.

Нехай $\Delta_{2n}=\{0=x_0<x_1<\dots<x_{2n}=2\pi\}$, $x_k=k\pi/n$ --- розбиття проміжку $[0,2\pi]$ та 
\begin{equation*}
\Psi_{\beta,1}(t)=(\Psi_{\beta}\ast B_{1})(t)=\sum\limits_{k=1}^{\infty}\frac{\psi(k)}{k}\cos\left(kt-\frac{(\beta+1)\pi}{2}\right),
\end{equation*}
де $\Psi_\beta(t)$ --- ядро вигляду \eqref{Psi_beta}, а $B_{1}(t)=\sum\limits_{k=1}^{\infty}k^{-1}\sin kt$ --- ядро Бернуллі.
Через $S\Psi_{\beta,1}(\Delta_{2n})$ позначатимемо простір $SK$-сплайнів $S\Psi_{\beta,1}(\cdot)$ за розбиттям $\Delta_{2n}$, тобто множину функцій виду
\begin{equation}\label{SK}
S\Psi_{\beta,1}(\cdot)=\alpha_{0}+\sum\limits_{k=1}^{2n} \alpha_k \Psi_{\beta,1}(\cdot-x_k),\;\sum\limits_{k=1}^{2n} \alpha_k=0,
\end{equation}
\begin{equation*}
\alpha_k\in\mathbb{R},\; k= 0, 1, \dots, 2n.
\end{equation*}
Фундаментальним $SK$-сплайном називають функцію $\overline{S\Psi}_{\beta,1}(\cdot)=\overline{S\Psi}_{\beta,1}(y, \cdot)$ виду \eqref{SK}, що задовольняє співвідношення
\begin{equation*}
\overline{S\Psi}_{\beta,1}(y, y_k)=\delta_{0,k}=
\begin{cases}
0, & k=\overline{1,2n-1},    \\
1, & k=0,
\end{cases}
\end{equation*}
де $y_k=x_k+y$, $x_k=k\pi/n$, $y\in[0,\frac{\pi}{n})$. 
Як було зазначено у роботі \cite{Serdyuk_1995}, серед $(\psi, \beta)$-похідних будь-якого сплайна виду \eqref{SK} існує функція, яка є сталою на кожному інтервалі $(x_k, x_{k+1})$. Саме таку функцію будемо розуміти під записом $(\overline{S\Psi}_{\beta,1}(\cdot))_\beta^\psi$.

\textbf{Означення.} \textit{Будемо казати, що для деякого дійсного числа $y$ і розбиття $\Delta_{2n}$ ядро $\Psi_\beta(\cdot)$ вигляду \eqref{Psi_beta} задовольняє умову $C_{y,2n}$ (і записувати $\Psi_\beta\in C_{y,2n}$), якщо для цього ядра існує єдиний фундаментальний сплайн $\overline{S\Psi}_{\beta,1}(y,\cdot)$ і для нього виконуються рівності}
\begin{equation*}
\textrm{sign}(\overline{S\Psi}_{\beta,1}(y,t_k))_\beta^\psi=(-1)^k\varepsilon e_k,\, k=\overline{0,2n-1},
\end{equation*}
\textit{де $t_k=(x_{k}+x_{k+1})/2,$ $e_k$ дорівнює або 0, або 1, а $\varepsilon$ приймає значення $\pm1$ і не залежить від $k$.}

\textbf{Теорема 5 (О.К.~Кушпель \cite{Kushpel_1988,Kushpel_1989}).} \textit{Нехай при деякому $n\in\mathbb{N}$ функція $\Psi_{\beta}$ вигляду \eqref{Psi_beta}, що породжує класи $C_{\beta,p}^\psi$, $p=1,\infty$, задовольняє умову $C_{y,2n}$, коли $y$ --- точка, в якій функція $|(\Psi_\beta\ast\varphi_n)(t)|$, $\varphi_n(t)=\textnormal{sign}\sin nt$, набуває максимального значення. Тоді }
\begin{equation*}
d_{2n}(C_{\beta,\infty}^\psi, C)\geqslant\|\Psi_\beta\ast\varphi_n\|_C,
\end{equation*}
\begin{equation*}
d_{2n-1}(C_{\beta,1}^\psi, L)\geqslant\|\Psi_\beta\ast\varphi_n\|_C.
\end{equation*}

Достатні умови включення $\Psi_\beta\in C_{y,2n}$ для ядер виду \eqref{Psi_beta} при тих чи інших обмеженнях на ядра $\Psi_\beta$ були встановлені у роботах \cite{Kushpel_1988,Serdyuk_1998,Serdyuk_1995,Serdyuk_1999}. Це дозволило авторам зазначених робіт застосувати теорему 5 і одержати в ряді нових випадків точні оцінки поперечників $d_m(C_{\beta,\infty}^\psi,C)$ та $d_m(C_{\beta,1}^\psi,L)$.

Для успішного застосування теореми 5 необхідно отримати певну інформацію про поведінку функцій $(\overline{S\Psi}_{\beta,1}(y,t))_\beta^\psi$. З цією метою встановимо наступне допоміжне твердження.

\textbf{Лема 1.} \textit{Нехай $\beta\in \mathbb{R}$, $\sum\limits_{k=1}^{\infty}\psi(k)<\infty$ і $y\in[0,\frac{\pi}{n})$ таке, що
\begin{equation}\label{lambda_not0}
|\lambda_{l}(y)|\not=0,\;l=\overline{1,n},
\end{equation}
де
\begin{equation}\label{lambd}
\lambda_{l}(y)=\frac1n \sum_{\nu=1}^{2n}e^{il\nu\pi/n}\Psi_{\beta,1}(y-\frac{\nu\pi}{n}).
\end{equation}
Тоді для довільного ${t\in(\frac{(k-1)\pi}{n},\frac{k\pi}{n})}$, $k=\overline{1,2n}$, виконується рівність}
\begin{equation*}
(\overline{S\Psi}_{\beta,1}(y,t))_\beta^\psi=
(-1)^{k+1}\frac{\pi}{4n\psi(n)}
\times
\end{equation*}
\begin{equation}\label{SP_Psi}
\times
\Bigg(
	\bigg(
	\frac{1}{2}
	+2\frac{\psi(n)}{n}\sum_{j=1}^{n-1}\frac{\cos j(t_k-y)}{|\lambda_{n-j}(y)|\cos\frac{j\pi}{2n}}
	\bigg)
	\mathop{\text{sign}} \sin(ny-\frac{\beta\pi}{2})
	+\gamma_1(y)
	+\gamma_2(y)
\Bigg),
\end{equation}
\textit{в якій $t_k=\frac{k \pi}{n}-\frac{\pi}{2n}$, а } 
\begin{equation}\label{gamma_2}
\gamma_1(y)
	=
		\gamma_1(\psi,\beta,k,y)
	=
		\frac{\psi(n)}{n}
		\left(
			\frac{z_{0}(y)}{|\lambda_{n}(y)|^2}+
			2\sum_{j=1}^{n-1}
				\frac{z_{j}(y)}{|\lambda_{n-j}(y)|^2\cos\frac{j\pi}{2n}}
		\right),
\end{equation}
\begin{equation}\label{gamma_3}
\gamma_2(y)
	=
		\gamma_2(\psi,\beta,y)
	=
		-
		\frac
			{R_0(y)\frac{n}{\psi(n)}}
			{2(2+R_0(y)\frac{n}{\psi(n)})}
		\mathop{\text{sign}} \sin(ny-\frac{\beta\pi}{2}),
\end{equation}
\begin{equation*}
z_{j}(y)
	=
		z_{j}(\psi,\beta,k,y)
	=
		|r_{j}(y)|\cos(j(t_k-y)+\arg(r_{j}(y)))
		-
\end{equation*}
\begin{equation}\label{z_nj}
		-
		R_{j}(y)\cos(j(t_k-y))
		\mathop{\text{sign}} \sin(ny-\frac{\beta\pi}{2}),
		\;j=\overline{0,n-1},
\end{equation}
\begin{equation}\label{|lambda_n-j|}
R_{j}(y)
	=
		R_{j}(\psi,\beta,y)
	=
		|\lambda_{n-j}(y)|
		-
		\frac{\psi(n-j)}{n-j}
		-
		\frac{\psi(n+j)}{n+j},
		\;j=\overline{0,n-1},
\end{equation}
\begin{align}
\label{r_n-j}
&r_{j}(y)
	=
		\sum_{\nu=1}^3
			r_{j}^{(\nu)}(y),
		\; j=\overline{0,n-1},
\\
\nonumber
&r_{j}^{(1)}(y)
	=
		r_{j}^{(1)}(\psi,\beta,y)
	=
		\frac
			{\psi(3n-j)e^{i(3ny-\frac{(\beta+1)\pi}{2})}}
			{3n-j}
		+
\\
\nonumber
&\phantom{r_{j}^{(1)}(y)=}
	+
	\sum\limits_{m=2}^{\infty}
		\left(
			\frac
				{\psi((2m+1)n-j)e^{i((2m+1)ny-\frac{(\beta+1)\pi}{2})}}
				{(2m+1)n-j}
			+
		\right.
\\
\label{r_nj_1}
&\phantom{r_{j}^{(1)}(y)=}
		\left.
			+
			\frac
				{\psi((2m-1)n+j)e^{-i((2m-1)ny-\frac{(\beta+1)\pi}{2})}}
				{(2m-1)n+j}
		\right),
\\
\label{r_nj_2}
&r_{j}^{(2)}(y)
	=
		r_{j}^{(2)}(\psi,\beta,y)
	=
		i
		\left(
			\frac{\psi(n+j)}{n+j}
			-
			\frac{\psi(n-j)}{n-j}
		\right)
		\cos(ny-\frac{\beta\pi}{2}),
\\
\nonumber
&r_{j}^{(3)}(y)
	=
		r_{j}^{(3)}(\psi,\beta,y)
	=
		\left(
			\frac{\psi(n-j)}{n-j}
			+
			\frac{\psi(n+j)}{n+j}
		\right)
		\times	
\\
\label{r_nj_3}
&\phantom{r_{j}^{(3)}(y)=}
	\times		
		(|\sin(ny-\frac{\beta\pi}{2})|-1)
		\mathop{\text{sign}} \sin(ny-\frac{\beta \pi}{2}).
\end{align}

\textbf{Доведення.} Будемо виходити із отриманого у роботі \cite{Serdyuk_1995} зображення функції $(\overline{S\Psi}_{\beta,1}(y,t))_\beta^\psi$, згідно з яким за умови $|\lambda_j(y)|\not=0$, $j=\overline{1, n}$, для довільного $t\in(x_{k-1},x_k)$ виконується рівність
\begin{equation}\label{SPsi_v0}
(\overline{S\Psi}_{\beta,1}(y,t))_\beta^\psi=\frac{\pi}{4n^2}\left(2\sum_{j=1}^{n-1}\frac{\sin jt_k\cdot\rho_j(y)-\cos jt_k\cdot\sigma_j(y)}{|\lambda_j(y)|^2\sin\frac{j\pi}{2n}}+\frac{(-1)^{k+1}\rho_n(y)}{|\lambda_n(y)|^2}\right),
\end{equation}
де
\begin{equation*}
\lambda_j(\cdot)=\frac1n \sum_{\nu=1}^{2n}e^{ij\nu\pi/n}\Psi_{\beta,1}(\cdot-\frac{\nu\pi}{n}),
\end{equation*}
$i$ --- уявна одиниця, $\rho_j(\cdot)=\mathop{\text{Re}}(\lambda_j(\cdot))$, $\sigma_j(\cdot)=\mathop{\text{Im}}(\lambda_j(\cdot))$, $t_k=\frac{k \pi}{n}-\frac{\pi}{2n}$.

Змінивши порядок підсумовування доданків у сумі в правій частині рівності \eqref{SPsi_v0}, маємо
\begin{equation*}
\sum_{j=1}^{n-1}\frac{\sin jt_k\cdot\rho_j(y)-\cos jt_k\cdot\sigma_j(y)}{|\lambda_j(y)|^2\sin\frac{j\pi}{2n}}=
\end{equation*}
\begin{equation*}
=\sum_{j=1}^{n-1}\frac{\sin (n-j)t_k\cdot\rho_{n-j}(y)-\cos (n-j)t_k\cdot\sigma_{n-j}(y)}{|\lambda_{n-j}(y)|^2\sin\frac{(n-j)\pi}{2n}}=
\end{equation*}
\begin{equation}\label{S_Psi}
=(-1)^{k+1}\sum_{j=1}^{n-1}\frac{\cos jt_k\cdot\rho_{n-j}(y)-\sin jt_k\cdot\sigma_{n-j}(y)}{|\lambda_{n-j}(y)|^2\cos\frac{j\pi}{2n}}.
\end{equation}

З урахуванням \eqref{SPsi_v0}  і \eqref{S_Psi} для фундаментального $SK$-сплайна $\overline{S\Psi}_{\beta,1}(y,t)$, 
 за умови $|\lambda_j(y)|\not=0$, $j=\overline{1, n}$, одержуємо зображення
\begin{equation*}
(\overline{S\Psi}_{\beta,1}(y,t))_\beta^\psi=
\end{equation*}
\begin{equation}\label{SP_v0}
=\frac{(-1)^{k+1}\pi}{4n^2}\left(2\sum_{j=1}^{n-1}\frac{\cos jt_k\cdot\rho_{n-j}(y)-\sin jt_k\cdot\sigma_{n-j}(y)}{|\lambda_{n-j}(y)|^2\cos\frac{j\pi}{2n}}+\frac{\rho_n(y)}{|\lambda_n(y)|^2}\right).      
\end{equation}
Покажемо, що величини $\lambda_{n-j}(y)$ виду \eqref{lambd} при $j=\overline{0,n-1}$ можна виразити наступним чином:
\begin{equation}\label{lambda_n-j}
\lambda_{n-j}(y)=
e^{-ijy}
\left(
	\left(
		\frac{\psi(n-j)}{n-j}+\frac{\psi(n+j)}{n+j}
	\right)
	\mathop{\text{sign}} \sin(ny-\frac{\beta \pi}{2})
	+r_{j}(y)
\right),
\end{equation}
де величини $r_{j}(y)$ задаються рівностями \eqref{r_n-j}.

Перепишемо ядро $\Psi_{\beta,1}$ у комплексній формі
\begin{equation*}
\Psi_{\beta,1}(t)
=(\Psi_{\beta}\ast B_1)(t)
=\sum_{k=1}^\infty\frac{\psi(k)}{k}\cos (kt-\frac{(\beta+1)\pi}{2})
= \frac 12 \sideset{}{'}\sum\limits_{k=-\infty}^{\infty}c_ke^{ikt},
\end{equation*}
де
\begin{equation}\label{c_k}
c_k=\frac{\psi(k)}{k}e^{-i\frac{(\beta+1)\pi}{2}},\;c_{-k}=\frac{\psi(k)}{k}e^{i\frac{(\beta+1)\pi}{2}},\;k\in\mathbb{N},
\end{equation}
а штрих біля знака суми означає, що при підсумовуванні відсутній доданок з нульовим номером.

Підставивши у \eqref{lambd} замість ядра $\Psi_{\beta,1}$ його розклад у комплексний ряд Фур’є, одержимо
\begin{equation*}
\lambda_l(y)
	=
		\frac1n \sum_{\nu=1}^{2n}e^{il\nu\pi/n}\frac 12 \sideset{}{'}\sum\limits_{k=-\infty}^{\infty}c_ke^{ik(y-\nu\pi/n)} 
	=
		\frac{1}{2n} \sum_{\nu=1}^{2n}\sideset{}{'}\sum\limits_{k=-\infty}^{\infty}c_k e^{i(ky+(l-k)\nu\pi/n)} 
	=
\end{equation*}
\begin{equation}\label{eq:lambda_dop1}
	=
		\frac{1}{2n} \sideset{}{'}\sum\limits_{k=-\infty}^{\infty}c_k e^{iky} \sum_{\nu=1}^{2n}e^{i((l-k)\nu\pi/n)}.
\end{equation}

Неважко переконатись, що
\begin{equation}\label{eq:lambda_dop2}
\sum_{\nu=1}^{2n}e^{i((l-k)\nu\pi/n)}=\begin{cases}
0,& \text{якщо } k\not=l-2mn, m\in\mathbb{Z};\\
2n,& \text{якщо } k=l-2mn, m\in\mathbb{Z}.
\end{cases}
\end{equation}

З \eqref{eq:lambda_dop1} та \eqref{eq:lambda_dop2} при $l=\overline{1,n}$ випливає наступне представлення:
\begin{equation*}
\lambda_l(y)=\sum\limits_{m=-\infty}^{+\infty}c_{l-2mn}e^{i(l-2mn)y}
=\sum\limits_{m=-\infty}^{+\infty}c_{2mn+l}e^{i(2mn+l)y}.
\end{equation*}

Звідси при $l=n-j\,$,   $j=\overline{0,n-1}$, отримуємо
\begin{equation*}
\lambda_{n-j}(y)=\sum\limits_{m=-\infty}^{+\infty}c_{(2m+1)n-j}e^{i((2m+1)n-j)y}=
\end{equation*}
\begin{equation}\label{lambda_n-j_y0}
=e^{-ijy}(c_{n-j}e^{iny}+c_{-(n+j)}e^{-iny}+r_{j}^{(1)}(y)).
\end{equation}

З урахуванням \eqref{c_k} перетворимо перші два доданки в \eqref{lambda_n-j_y0} наступним чином:
\begin{equation*}
c_{n-j}e^{iny}+c_{-(n+j)}e^{-iny}
	=
		\frac{\psi(n-j)}{n-j}e^{i(ny-\frac{(\beta+1)\pi}{2})}
		+\frac{\psi(n+j)}{n+j}e^{-i(ny-\frac{(\beta+1)\pi}{2})}=
\end{equation*}
\begin{equation*}
=\left(\frac{\psi(n-j)}{n-j}+\frac{\psi(n+j)}{n+j}\right)\cos(ny-\frac{(\beta+1)\pi}{2})+
\end{equation*}
\begin{equation*}
+i\left(\frac{\psi(n-j)}{n-j}-\frac{\psi(n+j)}{n+j}\right)\sin(ny-\frac{(\beta+1)\pi}{2})=
\end{equation*}
\begin{equation}\label{c_n-j-}
=\left(\frac{\psi(n-j)}{n-j}+\frac{\psi(n+j)}{n+j}\right)\sin(ny-\frac{\beta\pi}{2})+r_{j}^{(2)}(y).
\end{equation}

Записавши $\sin(ny-\frac{\beta\pi}{2})$ у вигляді
\begin{equation*}\label{sin_ny}
\sin\left(ny-\frac{\beta\pi}{2}\right)=
|\sin(ny-\frac{\beta\pi}{2})|\,
\mathop{\text{sign}} \sin(ny-\frac{\beta \pi}{2}),
\end{equation*}
з \eqref{c_n-j-} маємо
\begin{equation*}
c_{n-j}e^{iny}+c_{-(n+j)}e^{-iny}=
\left(
	\frac{\psi(n-j)}{n-j}+\frac{\psi(n+j)}{n+j}
\right)
\mathop{\text{sign}} \sin(ny-\frac{\beta\pi}{2})+
\end{equation*}
\begin{equation*}
+\left(
	\frac{\psi(n-j)}{n-j}+\frac{\psi(n+j)}{n+j}
\right)
(|\sin(ny-\frac{\beta\pi}{2})|-1)
\mathop{\text{sign}} \sin(ny-\frac{\beta \pi}{2})
+r_{j}^{(2)}(y)=
\end{equation*}
\begin{equation}\label{c_n-j}
=\left(
	\frac{\psi(n-j)}{n-j}
	+\frac{\psi(n+j)}{n+j}
\right)
\mathop{\text{sign}} \sin(ny-\frac{\beta\pi}{2})
+r_{j}^{(2)}(y)
+r_{j}^{(3)}(y).
\end{equation}

Рівності \eqref{lambda_n-j_y0} та \eqref{c_n-j} доводять формулу \eqref{lambda_n-j}.

Перетворимо чисельник кожного доданка в правій частині рівності (\ref{SP_v0}). Для цього, з урахуванням (\ref{lambda_n-j}), запишемо
\begin{equation*}
\rho_{n-j}(y)= \mathop{\text{Re}}(\lambda_{n-j}(y))=
\end{equation*}
\begin{equation}\label{rho_n-j}
=\left(\frac{\psi(n-j)}{n-j}+\frac{\psi(n+j)}{n+j}\right)\cos jy \;
\mathop{\text{sign}} \sin(ny-\frac{\beta\pi}{2})
+\mathop{\text{Re}}(e^{-ijy}r_{j}(y));
\end{equation}
\begin{equation*}
\sigma_{n-j}(y)=\mathop{\text{Im}}(\lambda_{n-j}(y))=
\end{equation*}
\begin{equation}\label{sigma_n-j}
=-\left(\frac{\psi(n-j)}{n-j}+\frac{\psi(n+j)}{n+j}\right)\sin jy \;
\mathop{\text{sign}} \sin(ny-\frac{\beta\pi}{2})
+\mathop{\text{Im}}(e^{-ijy}r_{j}(y)).
\end{equation}

Застосовуючи \eqref{rho_n-j} та \eqref{sigma_n-j}, отримуємо
\begin{equation*}
\cos jt_k\cdot\rho_{n-j}(y)
-\sin jt_k\cdot\sigma_{n-j}(y)=
\end{equation*}
\begin{equation*}
=
\left(
	\frac{\psi(n-j)}{n-j}
	+\frac{\psi(n+j)}{n+j}
\right)
\cos(j(t_k-y))
\,\mathop{\text{sign}} \sin(ny-\frac{\beta\pi}{2})+
\end{equation*}
\begin{equation*}
+\cos jt_k\cdot\mathop{\text{Re}} (e^{-ijy} r_{j}(y))
-\sin jt_k\cdot\mathop{\text{Im}} (e^{-ijy} r_{j}(y))=
\end{equation*}
\begin{equation*}
=
\left(
	\frac{\psi(n-j)}{n-j}
	+\frac{\psi(n+j)}{n+j}
	+R_{j}(y)
\right)
\cos(j(t_k-y)) 
\,\mathop{\text{sign}} \sin(ny-\frac{\beta\pi}{2})
+
\end{equation*}
\begin{equation}\label{cos+sin}
+z_{j}(y)=
|\lambda_{n-j}(y)|\cos(j(t_k-y)) 
\,\mathop{\text{sign}} \sin(ny-\frac{\beta\pi}{2})
+z_{j}(y),
\end{equation}
де
\begin{equation*}
z_{j}(y)=
\cos jt_k\cdot\mathop{\text{Re}}(e^{-ijy} r_{j}(y))
-\sin jt_k\cdot\mathop{\text{Im}} (e^{-ijy} r_{j}(y))-
\end{equation*}
\begin{equation*}
-R_{j}(y)\cos(j(t_k-y))
\,\mathop{\text{sign}} \sin(ny-\frac{\beta\pi}{2}),
\end{equation*}
а $R_{j}(y)$ означені у \eqref{|lambda_n-j|}.

В силу очевидної рівності
\begin{equation*}
e^{-ijy} r_{j}(y)=|r_{j}(y)|(\cos (\arg(r_{j}(y))-jy)+i\sin (\arg(r_{j}(y))-jy))
\end{equation*}
величину $z_{j}(y)$ можна зобразити у вигляді \eqref{z_nj}.

При $j=0$ формула \eqref{lambda_n-j} перетворюється в наступну рівність:
\begin{equation}\label{lambda_n}
\lambda_{n}(y)=
2\frac{\psi(n)}{n}
\mathop{\text{sign}} \sin(ny-\frac{\beta\pi}{2})
+r_{0}(y),
\end{equation}
де $r_{0}(y)$ визначається формулою \eqref{r_n-j}, у якій
\begin{equation}\label{r_n1}
r_{0}^{(1)}(y)=
2\sum\limits_{m=2}^{\infty}\frac{\psi((2m-1)n)}{(2m-1)n}\cos((2m-1)ny-\frac{(\beta+1)\pi}{2}),
\end{equation}
\begin{equation}\label{r_n2}
r_{0}^{(2)}(y)=0,
\end{equation}
\begin{equation}\label{r_n3}
r_{0}^{(3)}(y)=2\frac{\psi(n)}{n} (|\sin(ny-\frac{\beta\pi}{2})|-1)
\mathop{\text{sign}} \sin(ny-\frac{\beta\pi}{2}).
\end{equation}

З \eqref{lambda_n}--\eqref{r_n3} випливає, що $\sigma_{n}(y)=0$ і тому
\begin{equation*}
\rho_{n}(y)=\lambda_{n}(y)=
2\frac{\psi(n)}{n}
\mathop{\text{sign}} \sin(ny-\frac{\beta\pi}{2})
+r_{0}(y).  
\end{equation*}

Звідси, враховуючи \eqref{z_nj} та \eqref{|lambda_n-j|}, можна записати
\begin{equation*}
\rho_{n}(y)=
\left(
	2\frac{\psi(n)}{n}
	+R_0(y)
\right)
\mathop{\text{sign}} \sin(ny-\frac{\beta\pi}{2})
+z_{0}(y)=
\end{equation*}
\begin{equation}\label{rho_n}
=|\lambda_{n}(y)|
\mathop{\text{sign}} \sin(ny-\frac{\beta\pi}{2})
+z_{0}(y),
\end{equation}
де
\begin{equation*}
z_0(y)=r_0(y)- R_0(y)
\mathop{\text{sign}} \sin(ny-\frac{\beta\pi}{2}).
\end{equation*}

Із зображення (\ref{SP_v0}) і рівностей (\ref{cos+sin}) та \eqref{rho_n} отримуємо
\begin{equation*}
(\overline{S\Psi}_{\beta,1}(y,t))_\beta^\psi=
\end{equation*}
\begin{equation*}
=\frac{(-1)^{k+1}\pi}{4n\psi(n)}
\left(
	\mathop{\text{sign}} \sin(ny-\frac{\beta\pi}{2})
	\left(
		2\frac{\psi(n)}{n}\sum_{j=1}^{n-1}\frac{\cos j(t_k-y)}{|\lambda_{n-j}(y)|\cos\frac{j\pi}{2n}}
		+\frac{\psi(n)}{n|\lambda_n(y)|}
	\right)
	+
\right.
\end{equation*}
\begin{equation*}
\left.
	+2\frac{\psi(n)}{n}\sum_{j=1}^{n-1}\frac{z_{j}(y)}{|\lambda_{n-j}(y)|^2\cos\frac{j\pi}{2n}}
	+\frac{\psi(n)z_{0}(y)}{n|\lambda_{n}(y)|^2}
\right)=
\end{equation*}
\begin{equation*}
=\frac{(-1)^{k+1}\pi}{4n\psi(n)}
\Bigg(
	\mathop{\text{sign}} \sin(ny-\frac{\beta\pi}{2})
\Bigg.
\times
\end{equation*}
\begin{equation}\label{Phi_x}
\times
\Bigg.
	\Bigg(
		2\frac{\psi(n)}{n}\sum_{j=1}^{[\sqrt{n}]}\frac{\cos j(t_k-y)}{|\lambda_{n-j}(y)|\cos\frac{j\pi}{2n}}
		+\frac{\psi(n)}{n|\lambda_n(y)|}
	\Bigg)
	+\gamma_1(y)
\Bigg).
\end{equation}

В силу (\ref{|lambda_n-j|})
\begin{equation*}
\frac{\psi(n)\mathop{\text{sign}} \sin(ny-\frac{\beta\pi}{2})}{n|\lambda_n(y)|}=
\frac{\mathop{\text{sign}} \sin(ny-\frac{\beta\pi}{2})}{2+R_0(y)\frac{n}{\psi(n)}}=
\end{equation*}
\begin{equation*}
=\Bigg(
	\frac{1}{2}
	-\frac{R_0(y)\frac{n}{\psi(n)}}{2(2+R_0(y)\frac{n}{\psi(n)})}
\Bigg)
\mathop{\text{sign}} \sin(ny-\frac{\beta\pi}{2})=
\end{equation*}
\begin{equation}\label{s0}
=\frac{1}{2}\mathop{\text{sign}} \sin(ny-\frac{\beta\pi}{2})+\gamma_2(y).
\end{equation}

Із \eqref{Phi_x} та \eqref{s0} отримуємо \eqref{SP_Psi}.
Лему доведено.

Лема 1 дозволяє одержати зручне для подальших досліджень зображення величин $(\overline{S\Psi}_{\beta,1}(y,t))_\beta^\psi$, що породжуються ядрами $\Psi_\beta$ вигляду \eqref{Psi_beta}, коефіцієнти $\psi(k)$ яких належать множині $\mathcal{D}_q$, тобто виконуються умови \eqref{Dqumova}.

\textbf{Лема 2.} \textit{Нехай $\psi\in\mathcal{D}_q$, $q\in(0,1)$, $\beta\in \mathbb{R}$, $y\in[0,\frac{\pi}{n})$.
Тоді при виконанні умови \eqref{lambda_not0} для довільного ${t\in(\frac{(k-1)\pi}{n},\frac{k\pi}{n})}$, $k=\overline{1,2n}$, справедлива рівність}
\begin{equation}\label{SP_Phi_}
(\overline{S\Psi}_{\beta,1}(y,t))_\beta^\psi=
(-1)^{k+1}\frac{\pi}{4n\psi(n)}\;
\bigg(
	\mathcal{P}_q(t_k-y)
	\mathop{\text{sign}} \sin(ny-\frac{\beta\pi}{2})
	+\sum_{m=1}^5\gamma_m(y)
\bigg),
\end{equation}
\textit{в якій $t_k=\frac{k \pi}{n}-\frac{\pi}{2n}$, $\mathcal{P}_q(t)$ --- ядро аналітично продовжуваних в смугу функцій:}
\begin{equation*}
\mathcal{P}_q(t)= \frac{1}{2}+2\sum_{j=1}^{\infty}\frac{\cos jt}{q^j+q^{-j}},
\; q\in(0,1),
\end{equation*}
\textit{величини $\gamma_1(y)$ та $\gamma_2(y)$ задані рівностями \eqref{gamma_2} і \eqref{gamma_3} відповідно, а}
\begin{equation}\label{gamma_1}
\gamma_3(y)
	=
		\gamma_3(\psi,\beta,k,y)
	=
		2\sum_{j=[\sqrt{n}]+1}^{n-1}
			\frac
				{\cos j(t_k-y)}
				{\frac{n}{\psi(n)}|\lambda_{n-j}(y)|\cos\frac{j\pi}{2n}}
		\mathop{\text{sign}} \sin(ny-\frac{\beta\pi}{2}),
\end{equation}
\begin{equation}\label{gamma_4}
\gamma_4(y)
	=
		\gamma_4(\psi,\beta,k,y)
	=
		-2\sum_{j=1}^{[\sqrt{n}]}
			\frac
				{\delta_{j}(y)\cos j(t_k-y)}
				{\frac{n}{\psi(n)}|\lambda_{n-j}(y)|\cos\frac{j\pi}{2n}}
			\mathop{\text{sign}} \sin(ny-\frac{\beta\pi}{2}),
\end{equation}
\begin{equation}\label{gamma_5}
\gamma_5(y)
	=
		\gamma_5(q,\beta,k,y)
	=
		-2\sum\limits_{j=[\sqrt{n}]+1}^{\infty}
			\frac
				{\cos j(t_k-y)}
				{q^j+q^{-j}}
			\mathop{\text{sign}} \sin(ny-\frac{\beta\pi}{2}),
\end{equation}
\begin{equation}\label{delta_0}
\delta_{j}(y)
	=
		\delta_{j}(\psi,y)
	=
	\frac
		{n|\lambda_{n-j}(y)|\cos\frac{j\pi}{2n}}
		{(q^{-j}+q^{j})\psi(n)}
	-1,
	\;j=\overline{1,[\sqrt{n}]},
\end{equation}
\textit{$[a]$ --- ціла частина числа $a$.}

\textbf{Доведення.} Згідно з позначенням \eqref{gamma_1}
\begin{equation*}
2\frac{\psi(n)}{n}\sum_{j=1}^{n-1}
	\frac{\cos j(t_k-y)}{|\lambda_{n-j}(y)|\cos\frac{j\pi}{2n}}
	\mathop{\text{sign}} \sin(ny-\frac{\beta\pi}{2})
	=
\end{equation*}
\begin{equation}\label{summa_gamma}
=2\frac{\psi(n)}{n}\sum_{j=1}^{[\sqrt{n}]}\frac{\cos j(t_k-y)}{|\lambda_{n-j}(y)|\cos\frac{j\pi}{2n}}
\mathop{\text{sign}} \sin(ny-\frac{\beta\pi}{2})
+\gamma_3(y).
\end{equation}

Далі в силу формул \eqref{gamma_4}--\eqref{delta_0} можна записати рівності
\begin{equation*}
2\frac{\psi(n)}{n}\sum_{j=1}^{[\sqrt{n}]}\frac{\cos j(t_k-y)}{|\lambda_{n-j}(y)|\cos\frac{j\pi}{2n}}
\mathop{\text{sign}} \sin(ny-\frac{\beta\pi}{2})=
\end{equation*}
\begin{equation*}
=
2\sum_{j=1}^{[\sqrt{n}]}\frac{\cos j(t_k-y)}{(q^j+q^{-j})(1+\delta_{j}(y))}
\mathop{\text{sign}} \sin(ny-\frac{\beta\pi}{2})=
\end{equation*}
\begin{equation*}
=
\left(
	2\sum_{j=1}^{[\sqrt{n}]}
		\frac{\cos j(t_k-y)}{q^j+q^{-j}}
	-2\sum_{j=1}^{[\sqrt{n}]}
		\frac{\delta_{j}(y)\cos j(t_k-y)}{(q^j+q^{-j})(1+\delta_{j}(y))}
\right)
\mathop{\text{sign}} \sin(ny-\frac{\beta\pi}{2})
=
\end{equation*}
\begin{equation}\label{s1}
=\left(
	\mathcal{P}_q(t_k-y)-\frac{1}{2}
\right)
\mathop{\text{sign}} \sin(ny-\frac{\beta\pi}{2})	
+\gamma_4(y)+\gamma_5(y),
\end{equation}
Із \eqref{SP_Psi}, \eqref{summa_gamma} та \eqref{s1}  отримуємо \eqref{SP_Phi_}.
Лему доведено.

Як зазначалось вище, послідовності $\psi(k)=\frac{1}{\ch kh}$ ядра $H_{h,\beta}(t)$ вигляду \eqref{Psi_kernel} задовольняють умову $\mathcal{D}_q$ при $q=e^{-h}$, а тому для вказаних $\psi$ застосовна лема~2. Отже, при виконанні \eqref{lambda_not0} для $SK$-сплайнів, породжених ядром $H_{h,\beta}(t)$, має місце зображення \eqref{SP_Phi_}.

Наступне твердження містить оцінку зверху суми $\sum\limits_{k=1}^5|\gamma_k(y)|$ в зображенні \eqref{SP_Phi_} для ядер $\Psi_{\beta}(t)=H_{h,\beta}(t)$ у випадку, коли $y=y_0$, де $y_0$ --- точка, в якій фукнція $|\Psi_\beta\ast\varphi_n|$ набуває найбільшого значення.

\textbf{Лема 3.} \textit{Нехай величини $\gamma_l(y_0)$, $l=\overline{1,5}$, задаються рівностями \eqref{gamma_2}, \eqref{gamma_3}, \eqref{gamma_1}--\eqref{gamma_5}, в яких $\psi(n)=\frac{1}{\ch nh}=\frac{2q^n}{1+q^{2n}}$, $h>0$, $q=e^{-h}$, $\beta\in\mathbb{R}$, а $y_0=y_0(n,h,\beta)=\frac{\theta_n\pi}{n}$, де $\theta_n$ --- корінь рівняння \eqref{theta}, $\theta_n\in[0,1)$. Тоді при $n\geqslant 9$ та виконанні умови }
\begin{equation}\label{umova_z}
\frac{q^{n}}{1-q^{2n}}
	\leqslant 
			\frac {7q^{\sqrt{n}}}{37n^2}
\end{equation}
\textit{для довільного $t\in(\frac{(k-1)\pi}{n},\frac{k\pi}{n})$, $k=\overline{1,2n}$, має місце зображення }
\begin{equation}\label{SP_Phi_q}
(\overline{S\Psi}_{\beta,1}(y_0,t))_\beta^\psi=
(-1)^{k+1}\frac{\pi}{4n\psi(n)}\;
\bigg(
	\mathcal{P}_q(t_k-y_0)
	\mathop{\text{sign}} \sin(ny_0-\frac{\beta\pi}{2})
	+\sum_{l=1}^5\gamma_l(y_0)
\bigg),
\end{equation}
\textit{ та справедлива оцінка}
\begin{equation*}
\sum\limits_{l=1}^5
	|\gamma_l(y_0)|
	\leqslant
		\frac{37}{5(1-q)}q^{\sqrt{n}}
		+
		\frac{q}{(1-q)^2}
		\min
			\left\{
				\frac{160}{27(n-\sqrt{n})},
				\frac{8}{3n-7\sqrt{n}}
			\right\}.
\end{equation*}

\textbf{Доведення.} Встановимо спочатку справедливість зображення \eqref{SP_Phi_q}. Для цього достатньо показати, що при $y=y_0$ виконується умова \eqref{lambda_not0}. Знайдемо оцінки зверху величин $|r_{j}(y_0)|$ та $|R_{j}(y_0)|$ при $j=\overline{0,n-1}$.
З \eqref{r_nj_1} маємо
\begin{equation*}
|r_{j}^{(1)}(y_0)|
	\leqslant
		\frac{2q^{3n-j}}{(3n-j)(1+q^{2(3n-j)})}
		+
		2\sum\limits_{m=2}^{\infty}
			\left(
				\frac{q^{(2m+1)n-j}}{((2m+1)n-j)(1+q^{2((2m+1)n-j)})}
				+				
			\right.
\end{equation*}
\begin{equation*}
			\left.
				+
				\frac{q^{(2m-1)n+j}}{((2m-1)n+j)(1+q^{2((2m-1)n+j)})}
			\right)
	<
\end{equation*}
\begin{equation*}
	<
		\frac{2q^{3n-j}}{3n-j}
		+
		2\sum\limits_{m=2}^{\infty}
			\left(
				\frac{q^{(2m+1)n-j}}{(2m+1)n-j}
				+
				\frac{q^{(2m-1)n+j}}{(2m-1)n+j}
		\right)
	=
\end{equation*}
\begin{equation}\label{r1_addit}
	=
		2\sum\limits_{m=1}^{\infty}
			\left(
				\frac{q^{(2m+1)n-j}}{(2m+1)n-j}
				+
				\frac{q^{(2m+1)n+j}}{(2m+1)n+j}
			\right).
\end{equation}
Враховуючи, що внаслідок опуклості послідовності $\frac{q^k}{k}$ виконується нерівність ${\frac{q^{k-j}}{k-j}+\frac{q^{k+j}}{k+j}<\frac{q^{k-n}}{k-n}+\frac{q^{k+n}}{k+n}}$, $k>n$, $j=\overline{0,n-1}$, із \eqref{r1_addit} знаходимо
\begin{equation*}
|r_{j}^{(1)}(y_0)|
	\leqslant
		2\sum\limits_{m=1}^{\infty}
			\left(
				\frac{q^{2mn}}{2mn}
				+
				\frac{q^{2(m+1)n}}{2(m+1)n}
			\right)
	=
\end{equation*}
\begin{equation}\label{|r_nj1|}
	= 
		\frac{q^{2n}}{n}
		+
		2\sum\limits_{m=2}^{\infty}
			\frac{q^{2mn}}{mn}
	\leqslant
		\frac{1}{n}
		\sum\limits_{m=1}^{\infty}
			q^{2mn}
	=
		\frac{q^{2n}}{n(1-q^{2n})}.
\end{equation}

З рівняння \eqref{theta} при $q=e^{-h}$ отримуємо
\begin{equation*}
\left|
	\cos
	\left(
		\theta_n\pi-\frac{\beta\pi}{2}
	\right)
\right|
	=
		(1+q^{2n})
		\left|
			\sum\limits_{\nu=1}^{\infty}
				\frac{q^{2\nu n}}{1+q^{2(2\nu+1) n}}
				\cos
				\left(
					(2\nu+1)\theta_n\pi-\frac{\beta\pi}{2}
				\right)
		\right|
	<
\end{equation*}
\begin{equation}\label{cos_y}
	< 
		(1+q^{2n})
		\sum\limits_{\nu=1}^{\infty}
			q^{2\nu n}
	=
		\frac
			{q^{2n}}
			{1-q^{2n}}
		(1+q^{2n}).
\end{equation}

У свою чергу, з \eqref{cos_y} випливає, що
\begin{equation}\label{a_n}
0
	\leqslant
		1
		-
		|\sin(\theta_n\pi-\frac{\beta\pi}{2})|
	\leqslant 
		|\cos(\theta_n\pi-\frac{\beta\pi}{2})|
	< 
		\frac
			{q^{2n}}
			{1-q^{2n}}
		(1+q^{2n}).
\end{equation}

З \eqref{cos_y} та  \eqref{r_nj_2} маємо
\begin{equation}\label{|r_nj2|}
|r_{j}^{(2)}(y_0)|
	<
		\frac
			{2q^{2n}(1+q^{2n})}
			{1-q^{2n}}
		\left(
			\frac
				{q^{n-j}}
				{(n-j)(1+q^{2(n-j)})}
			-
			\frac
				{q^{n+j}}
				{(n+j)(1+q^{2(n+j)})}
		\right).
\end{equation}

Із \eqref{a_n} та \eqref{r_nj_3} знаходимо
\begin{equation}\label{|alpha_n|}
|r_{j}^{(3)}(y_0)|
	<
		\frac
			{2q^{2n}(1+q^{2n})}
			{1-q^{2n}}
		\left(
			\frac
				{q^{n-j}}
				{(n-j)(1+q^{2(n-j)})}
			+
			\frac
				{q^{n+j}}
				{(n+j)(1+q^{2(n+j)})}
		\right).
\end{equation}

З умови \eqref{umova_z} випливає, що
\begin{equation}\label{q2n}
q^{2n}<\frac{49}{1369n^4}.       
\end{equation}
Отже, з \eqref{|r_nj1|}, \eqref{|r_nj2|}, \eqref{|alpha_n|} та \eqref{q2n} при $n\geqslant 9$ випливає оцінка величини $|r_{j}(y_0)|$ 
\begin{equation*}
|r_{j}(y_0)|
	\leqslant 
		|\sum\limits_{\nu=1}^{3}
			r_{j}^{(\nu)}(y_0)|
	< 
		\frac
			{q^{2n}}
			{1-q^{2n}}
		\left(
			\frac
				{4(1+q^{2n})q^{n-j}}
				{(n-j)(1+q^{2(n-j)})}			
			+
			\frac{1}{n}
		\right)
	<
\end{equation*}
\begin{equation}\label{|r_nj|}
	< 
		\frac
			{q^{2n}}
			{1-q^{2n}}
		\left(
			4q(1+q^{2n})
			+
			\frac{1}{n}
		\right)
	< 
		\frac{38}{9}
		\frac
			{q^{2n}}
			{1-q^{2n}},
j=\overline{0,n-1}.
\end{equation}

При $j=0$ оцінку \eqref{|r_nj|} можна покращити. Дійсно, в силу \eqref{r_n1} маємо
\begin{equation*}
|r_{0}^{(1)}(y_0)|
	\leqslant
		4\sum\limits_{m=2}^{\infty}
			\frac
				{q^{(2m-1)n}}
				{((2m-1)n)(1+q^{2(2m-1)n})}
	<
		\frac{4}{3n}
		\sum\limits_{m=2}^{\infty}
			q^{(2m-1)n}
	=
		\frac{4}{3n}
		\frac
			{q^{3n}}
			{1-q^{2n}},
\end{equation*}
а з \eqref{r_n3} та \eqref{a_n} випливає
\begin{equation*}
\left| 
	r_{0}^{(3)}(y_0)
\right|
	< 
		\frac
			{4q^{3n}}
			{n(1-q^{2n})}.
\end{equation*}
Тоді, враховуючи \eqref{r_n2}, можемо записати
\begin{equation}\label{|r_n|}
|r_{0}(y_0)|
	\leqslant 
		|r_{0}^{(1)}(y_0)
		+
		r_{0}^{(3)}(y_0)|
	\leqslant 
		\frac{16}{3n} 
		\frac
			{q^{3n}}
			{1-q^{2n}}.
\end{equation}

Із \eqref{lambda_n-j} отримуємо зображення
\begin{equation*}
|\lambda_{n-j}(y_0)|
	=		
\end{equation*}
\begin{equation*}
	=
		\left|
			\mathop{\text{sign}} \sin(ny-\frac{\beta\pi}{2})
			\left(
				\frac
					{2q^{n-j}}
					{(n-j)(1+q^{2(n-j)})}
				+
				\frac
					{2q^{n+j}}
					{(n+j)(1+q^{2(n+j)})}
			\right)
			+
			r_{j}(y_0)
		\right|,
\end{equation*}
з якого безпосередньо випливає оцінка
\begin{equation}\label{|lambda_n-j|0'}
|\lambda_{n-j}(y_0)|
	\leqslant 
		\frac
			{2q^{n-j}}
			{(n-j)(1+q^{2(n-j)})}
		+
		\frac
			{2q^{n+j}}
			{(n+j)(1+q^{2(n+j)})}
		+
		|r_{j}(y_0)|.
\end{equation}
Оскільки внаслідок \eqref{a_n} та умови \eqref{umova_z}
\begin{equation}\label{sing_not_zero}
|\sin(ny_0-\frac{\beta\pi}{2})|
	\geqslant
		1
		-
		\frac
			{q^{2n}}
			{1-q^{2n}}
		(1+q^{2n})
	>0,
\end{equation}
то отримуємо також оцінку
\begin{equation}\label{|lambda_n-j|0''}
|\lambda_{n-j}(y_0)|
	\geqslant 
		\frac
			{2q^{n-j}}
			{(n-j)(1+q^{2(n-j)})}
		+
		\frac
			{2q^{n+j}}
			{(n+j)(1+q^{2(n+j)})}
		-
		|r_{j}(y_0)|.
\end{equation}

В силу \eqref{|lambda_n-j|}, \eqref{|lambda_n-j|0'} та \eqref{|lambda_n-j|0''}
\begin{equation}\label{|R_n,j|}
|R_{j}(y_0)|
	\leqslant 
		|r_{j}(y_0)|,
\; j=\overline{0,n-1}.
\end{equation}

Взявши до уваги оцінки  \eqref{|lambda_n-j|0''} та \eqref{|r_nj|}, маємо
\begin{equation*}
|\lambda_{n-j}(y_0)|
	>
		\frac
			{q^{n-j}}
			{n-j}
		+
		\frac
			{q^{n+j}}
			{n+j}
		-
		\frac{38}{9}
		\frac
			{q^{2n}}
			{1-q^{2n}}		
	=
\end{equation*}
\begin{equation}\label{|lambda_n-j|3}
	=
		\frac
			{q^{n}}
			{n-j} 
		\left(
			q^{-j}
			+
			\frac{n-j}{n+j}
			q^{j}
			-
			\frac{38(n-j)q^{n}}{9(1-q^{2n})}
		\right).
\end{equation}

Оскільки при $j=\overline{0,n-1}$ та $n\geqslant 9\;\;$ 
$\dfrac{9q^{-j}}{380(n-j)}
>\dfrac{9}{380n}
>\dfrac{7q^{\sqrt{n}}}{37n^2}$, то з умови \eqref{umova_z} випливає нерівність
\begin{equation*}
\frac{9}{380(n-j)}
q^{-j}
	>
		\frac{q^{n}}{1-q^{2n}},
\end{equation*}
яка еквівалентна наступній нерівності:
\begin{equation}\label{um1}
\frac
	{q^{-j}}
	{10}
	>
		\frac{38(n-j)q^{n}}{9(1-q^{2n})},
\; j=\overline{0,n-1}.
\end{equation}

В силу \eqref{um1} виконуються оцінки
\begin{equation*}
q^{-j}
+
\frac{n-j}{n+j}
q^{j}
-
\frac{38(n-j)q^{n}}{9(1-q^{2n})}
	=
\end{equation*}
\begin{equation}\label{ots}
	=
		\frac{9q^{-j}}{10}
		+
		\frac{q^{-j}}{10}
		+
		\frac{n-j}{n+j}
		q^{j}
		-
		\frac{38(n-j)q^{n}}{9(1-q^{2n})}
	>
		\frac{9q^{-j}}{10},
\; j=\overline{0,n-1}.
\end{equation}

Об’єднуючи \eqref{|lambda_n-j|3} та \eqref{ots}, маємо
\begin{equation}\label{mod_lambda_nj}
|\lambda_{n-j}(y_0)|
	>
		\frac
			{9q^{n-j}}
			{10(n-j)}.
\end{equation}

З нерівності \eqref{mod_lambda_nj} випливає виконання умови \eqref{lambda_not0}, а отже, і справедливість зображення \eqref{SP_Phi_q}. 

Знайдемо оцінки зверху кожної з величин $|\gamma_l(y_0)|$, $l=\overline{1,5}$. 
Розпочнемо з оцінки величини $|\gamma_1(y_0)|$.
Оскільки для $x\in[0,\frac \pi2)$ справджується нерівність $\cos x\geqslant 1-\frac {2x}{\pi}>0$, отримуємо співвідношення
\begin{equation}\label{cos1}
\cos \frac{j\pi}{2n}
	\geqslant 
		1
		-
		\frac {j}{n}
	=
		\frac {n-j}{n},
\; j=\overline{0,n-1}.
\end{equation}

З \eqref{mod_lambda_nj} та \eqref{cos1} маємо
\begin{equation}\label{|lambda_n-j|_2}
\frac
	{n(1+q^{2n})}
	{2q^{n}}
|\lambda_{n-j}(y_0)|^2
\cos\frac{j\pi}{2n}
	> 
		\frac{81(1+q^{2n})}{200n}q^{n-2j}.
\end{equation}

З \eqref{|R_n,j|} і \eqref{z_nj} випливає, що $|z_{j}(y_0)|\leqslant 2|r_{j}(y_0)|$. Тому враховуючи \eqref{|r_nj|}, \eqref{|lambda_n-j|_2} та умову \eqref{umova_z}, з \eqref{gamma_2} одержуємо
\begin{equation*}
\left|\gamma_1(y_0)\right|
	\leqslant
		\frac{800}{81(1+q^{2n})}
		\max_{0\leqslant j\leqslant n-1}| r_{j}(y_0)| 
		\frac{n}{q^{n}}
		\sum_{j=0}^{n-1}
			q^{2j}
	<
\end{equation*}
\begin{equation}\label{r1}
	<
		\frac
			{30400\,nq^{n}}
			{729(1+q^{2n})(1-q^{2n})}
		\sum_{j=0}^{\infty}
			q^{2j}
	\leqslant
		\frac{212800}{26973n}
		q^{\sqrt{n}}
		\;\frac{1}{1-q^2}.
\end{equation}

Оцінимо $|\gamma_2(y_0)|$. З \eqref{gamma_3}, \eqref{|R_n,j|}, \eqref{|r_n|}, \eqref{umova_z} і \eqref{q2n} отримуємо
\begin{equation*}
|\gamma_2(y_0)|
	\leqslant
		\frac
			{\frac
				{8q^{2n}(1+q^{2n})}
				{3(1-q^{2n})}
			}
			{2\left|2
			-
			\frac
				{8q^{2n}(1+q^{2n})}
				{3(1-q^{2n})}\right|}
	=
		\frac
			{2q^{2n}(1+q^{2n})}
			{\left|
				3
				-
				7q^{2n}
				-
				4q^{4n}
			\right|}
	<
		\frac
			{4q^{2n}}
			{3-11q^{2n}}	
	=
\end{equation*}
\begin{equation*}
	=
		\frac
			{1-q^{2n}}
			{3-11q^{2n}}
		\;\frac{4q^{2n}}{1-q^{2n}}
	=
		\left(
			\frac {1}{11}
			+
			\frac{8}{11(3-11q^{2n})}
		\right)
		\frac{4q^{2n}}{1-q^{2n}}
	<
\end{equation*}
\begin{equation}\label{r2_2}
	<
		\frac {5}{11} 
		\frac{4q^{2n}}{1-q^{2n}}
	<
		\frac
			{140q^{n+\sqrt{n}}}
			{407n^2}.
\end{equation}

Оцінимо величину $|\gamma_3(y_0)|$.
Взявши до уваги \eqref{mod_lambda_nj} та \eqref{cos1}, маємо
\begin{equation}\label{|lambda_n-j|2}
\frac
	{n(1+q^{2n})}
	{2q^{n}}
|\lambda_{n-j}(y_0)|
\cos\frac{j\pi}{2n}
	> 
		\frac
			{9(1+q^{2n})q^{-j}}
			{20}.
\end{equation}

Тому, зважаючи на \eqref{|lambda_n-j|2}, з \eqref{gamma_1} знаходимо
\begin{equation}\label{r3}
\left|\gamma_3(y_0)\right|
	<
		\frac{40}{9(1+q^{2n})}
		\sum_{j=[\sqrt{n}]+1}^{n-1}
			q^j
	=
		\frac
			{40(q^{[\sqrt{n}]+1}-q^{n})}
			{9(1-q)(1+q^{2n})}
	\leqslant
		\frac
			{40q^{\sqrt{n}}}
			{9(1-q)}.
\end{equation}

Перш ніж оцінити $|\gamma_4(y_0)|$ встановимо оцінки зверху для величини $|\delta_j(y_0)|$, означеної в \eqref{delta_0}. З урахуванням \eqref{|lambda_n-j|}
\begin{equation*}
\frac{n(1+q^{2n})}{2q^{n}}
	|\lambda_{n-j}(y_0)|
		\cos\frac{j\pi}{2n}
	=
\end{equation*}
\begin{equation*}
	=		
		\left(
			\frac{n}{n-j}
			\frac
				{(1+q^{2n})q^{n-j}}
				{q^{n}(1+q^{2(n-j)})}
			+
			\frac{n}{n+j}
			\frac
				{(1+q^{2n})q^{n+j}}
				{q^{n}(1+q^{2(n+j)})}
			+
			R_{j}(y_0)
			\frac{(1+q^{2n})n}{2q^{n}}
		\right)
		\cos\frac{j\pi}{2n}
	=
\end{equation*}
\begin{equation*}
	=
		(q^{-j}+q^{j})
		(1-2\sin^2\frac{j\pi}{4n})
	+
		\left(
			\frac{n}{n-j}
			\frac
				{q^{-j}(1+q^{2n})}
				{1+q^{2(n-j)}}
			-
			q^{-j}
			+
		\right.
\end{equation*}
\begin{equation}\label{cos_lambda0}
		\left.
			+
			\frac{n}{n+j}
			\frac
				{q^{j}(1+q^{2n})}
				{1+q^{2(n+j)}}
			-
			q^{j}
			+
			R_{j}(y_0)
			\frac{n(1+q^{2n})}{2q^{n}}
		\right)
		\cos\frac{j\pi}{2n}.
\end{equation}

Оскільки
\begin{equation*}
\left|
	\frac{n}{n-j}
	\frac
		{1+q^{2n}}
		{1+q^{2(n-j)}}
	-
	1
\right|
	=
	\left|
		\left(
			1+\frac{j}{n-j}
		\right)
		\left(
		1-
		\frac
			{q^{2(n-j)}-q^{2n}}
			{1+q^{2(n-j)}}
		\right)
		-
		1
	\right|
	=
\end{equation*}
\begin{equation*}
	=
	\left|
		\frac{j}{n-j}
		-
		\frac{n}{n-j}
		\frac
			{q^{2(n-j)}-q^{2n}}
			{1+q^{2(n-j)}}
	\right|
	<
		\frac{j}{n-j}
		+
		\frac{n}{n-j}
		q^{2(n-j)},
\end{equation*}
\begin{equation*}
\left|
	\frac{n}{n+j}
	\frac
		{1+q^{2n}}
		{1+q^{2(n+j)}}
	-
	1
\right|
	=
	\left|
		\left(
			1-\frac{j}{n+j}
		\right)
		\left(
		1+
		\frac
			{q^{2n}-q^{2(n+j)}}
			{1+q^{2(n+j)}}
		\right)
		-
		1
	\right|
	=
\end{equation*}
\begin{equation*}
	=
	\left|
		-\frac{j}{n+j}
		+
		\frac{n}{n+j}
		\frac
			{q^{2n}-q^{2(n+j)}}
			{1+q^{2(n+j)}}
	\right|
	<
		\frac{j}{n+j}
		+
		\frac{n}{n+j}
		q^{2n},
\end{equation*}
то
\begin{equation}\label{max_frac}
\max\{
	|\frac{n}{n-j}
	\frac
		{1+q^{2n}}
		{1+q^{2(n-j)}}
	-
	1|,
	|\frac{n}{n+j}
	\frac
		{1+q^{2n}}
		{1+q^{2(n+j)}}
	-
	1|
\}
	<
		\frac{j}{n-j}
		+
		\frac{n}{n-j}
		q^{2(n-j)}.
\end{equation}
Тоді з \eqref{delta_0}, \eqref{|r_nj|}, \eqref{|R_n,j|}, \eqref{cos_lambda0}, \eqref{max_frac} та  з урахуванням опуклості послідовності $q^k$ для величин $|\delta_{j}(y_0)|$ будемо мати
\begin{equation*}
|\delta_{j}(y_0)|
	\leqslant 
		2\sin^2\frac{j\pi}{4n}
		+
		\frac{1}{q^{-j}+q^{j}}
		\left(
			\bigg(
				\frac{j}{n-j}
				+
				\frac
					{nq^{2(n-j)}}
					{n-j}
			\bigg)
			(q^{-j}
			+
			q^{j})
			+
		\right.
\end{equation*}
\begin{equation*}
		\left.
			+
			|R_{j}(y_0)|
			\frac
				{n(1+q^{2n})}
				{2q^{n}}
		\right)
	\leqslant 
		2\left(
			\frac{j\pi}{4n}
		\right)^2
		+
		\frac{j}{n-j}
		+
		\frac
			{nq^{2(n-j)}}
			{n-j}
		+
		\frac
			{n(1+q^{2n})|r_{j}(y_0)|}
			{2(q^{n-j}+q^{n+j})}
	\leqslant
\end{equation*}
\begin{equation*}
	\leqslant 
		\frac{j^2\pi^2}{8n^2}
		+
		\frac{j}{n-j}
		+
		\frac
			{nq^{2(n-j)}}
			{n-j}
		+
		\frac{19n(1+q^{2n})}{18}
		\frac
			{q^{n}}
			{1-q^{2n}}
	=
\end{equation*}
\begin{equation}\label{delta0}
	= 
		\frac{4j}{3(n-j)}
		+
		\left(
			\frac{j^2\pi^2}{8n^2}
			+
			\frac
				{nq^{2(n-j)}}
				{n-j}
			+
			\frac{19n(1+q^{2n})}{18}
			\frac
				{q^{n}}
				{1-q^{2n}}
			-
			\frac{j}{3(n-j)}
		\right).
\end{equation}

Покажемо, що  для усіх $j=\overline{1,[\sqrt{n}]}$
\begin{equation}\label{delta}
|\delta_{j}(y_0)|\leqslant \frac{4j}{3(n-j)}.
\end{equation}

Для цього, в силу  \eqref{delta0} досить переконатися, що при $j=\overline{1,[\sqrt{n}]}$ має місце нерівність
\begin{equation}\label{umova4}
\frac{j}{3(n-j)}
-
\frac{j^2\pi^2}{8n^2}
-
\frac
	{\sqrt{n}q^{2(n-\sqrt{n})}}
	{\sqrt{n}-1}
	>
		\frac{19n(1+q^{2n})}{18}
		\frac
			{q^{n}}
			{1-q^{2n}}.
\end{equation}

Дійсно, як показано у роботі \cite[с.~104]{My_JAT}, при кожному фіксованому ${x\geqslant9}$ функція $f(x,\tau)=\frac{\tau}{3(x-\tau)}-\frac{\tau^2\pi^2}{8x^2}$ на $[1, \sqrt{x}]$ набуває найменшого значення у точці $\tau=1$. Тому при $n\geqslant9$, з урахуванням \eqref{umova_z}, маємо для всіх $j=\overline{1,[\sqrt{n}]}$
\begin{equation*}
\frac{j}{3(n-j)}
-
\frac{j^2\pi^2}{8n^2}
-
\frac
	{\sqrt{n}q^{2(n-\sqrt{n})}}
	{\sqrt{n}-1}
	\geqslant
		\frac{1}{3(n-1)}
		-
		\frac{\pi^2}{8n^2}
		-
		\frac
			{\sqrt{n}q^{2(n-\sqrt{n})}}
			{2}
	>
\end{equation*}
\begin{equation}\label{dentaest}
	>
		\frac{1}{3(n-1)}
		-
		\frac{\pi^2}{8n^2}
		-
		\frac{49}{2738n^3\sqrt{n}} .
\end{equation}

При $n=9$, враховуючи \eqref{umova_z} та \eqref{q2n}, отримуємо
\begin{equation*}
\frac{1}{3(n-1)}
		-
\frac{\pi^2}{8n^2}
-
\frac{49}{2738n^3\sqrt{n}}
	=
		\frac{1}{24}
		-
		\frac{\pi^2}{648}
		-
		\frac{49}{5988006}           
	>
		\frac{17}{648}
		-
		\frac{49}{5988006}
	>
\end{equation*}
\begin{equation}\label{dentaestadd1}
	>
		0{,}025
	>
		\frac{7\cdot 19(1+q^{18})q^{3}}{18\cdot 37\cdot 9}
	=
		\frac{7\cdot 19(1+q^{2n})q^{\sqrt{n}}}{18\cdot 37n}
	>
		\frac{19n(1+q^{2n})q^{n}}{18(1-q^{2n})}.
\end{equation}

При $n\geqslant10$, враховуючи \eqref{umova_z} та \eqref{q2n}, отримуємо
\begin{equation*}
\frac{1}{3(n-1)}
		-
\frac{\pi^2}{8n^2}
-
\frac{49}{2738n^3\sqrt{n}}
	>
		\frac{1}{n}
		\left(
			\frac{1}{3}-\frac{\pi^2}{80}
		\right)
		-
		\frac{49}{2738n^3\sqrt{n}}           
	>
\end{equation*}
\begin{equation}\label{dentaestadd2}      
	>
		\frac{5}{24n}
		-
		\frac{49}{2738n^3\sqrt{n}} 
	>
		\frac{7\cdot 19(1+q^{2n})q^{\sqrt{n}}}{18\cdot 37n}
	>
		\frac{19n(1+q^{2n})q^{n}}{18(1-q^{2n})}.
\end{equation}

З \eqref{dentaest}, \eqref{dentaestadd1} та \eqref{dentaestadd2} випливає справедливість \eqref{umova4}, а отже, і \eqref{delta}.

Формули \eqref{gamma_4}, \eqref{|lambda_n-j|2} та \eqref{delta} дозволяють одержати 
при $n\geqslant9$ наступну оцінку величини $\gamma_4(y_0)$:
\begin{equation*}
|\gamma_4(y_0)|
	\leqslant
		2\sum_{j=1}^{[\sqrt{n}]}
			\frac
				{\frac{4j}{3(n-j)}}
				{\frac
					{9(1+q^{2n})q^{-j}}
					{20}}
	=
		\frac{160}{27(1+q^{2n})}
		\sum_{j=1}^{[\sqrt{n}]}
			\frac{j}{n-j}q^j
	\leqslant
\end{equation*}
\begin{equation}\label{|R|}
	\leqslant
		\frac{160}{27(n-\sqrt{n})}
		\sum_{j=1}^{[\sqrt{n}]}
			jq^j
	<
		\frac{160}{27(n-\sqrt{n})}
		\sum_{j=1}^{\infty}
			jq^j
	<
		\frac{160}{27(n-\sqrt{n})}
		\; \frac{q}{(1-q)^2}.
\end{equation}

Водночас для величини $|\gamma_4(y_0)|$ можна отримати іншу оцінку зверху. З цією метою, помітивши, що в силу \eqref{delta_0}
\begin{equation*}
\frac{n(1+q^{2n})}{2q^{n}}
|\lambda_{n-j}(y_0)|
\cos\frac{j\pi}{2n}
	=
		(q^j+q^{-j})
		(1+\delta_{j}(y_0)),
\end{equation*}
з \eqref{gamma_4} та \eqref{delta} при $n\geqslant9$ одержуємо
\begin{equation*}
|\gamma_4(y_0)|
	\leqslant
		2\sum_{j=1}^{[\sqrt{n}]}
			\frac
				{\frac{4j}{3(n-j)}}
				{1-\frac{4j}{3(n-j)}}
			q^j
	=
		2\sum_{j=1}^{[\sqrt{n}]}
			\frac{4j}{3n-7j}
			q^j
	\leqslant
\end{equation*}
\begin{equation}\label{|R|_}
	\leqslant
		\frac{8}{3n-7\sqrt{n}}
		\sum_{j=1}^{[\sqrt{n}]}
			jq^j
	<
		\frac{8}{3n-7\sqrt{n}}
		\sum_{j=1}^{\infty}
			jq^j
	<
		\frac{8}{3n-7\sqrt{n}}
		\; \frac{q}{(1-q)^2}.
\end{equation}

Із \eqref{|R|} і \eqref{|R|_} випливає оцінка
\begin{equation}\label{min|R|}
|\gamma_4(y_0)|
	\leqslant
		\frac{q}{(1-q)^2}
		\min
			\left\{
				\frac{160}{27(n-\sqrt{n})},
				\frac{8}{3n-7\sqrt{n}}
			\right\}.
\end{equation}

В силу \eqref{gamma_5} для величини $|\gamma_5(y_0)|$ маємо
\begin{equation}\label{r4}
\left|\gamma_5(y_0)\right|
	\leqslant 
		2\sum\limits_{j=[\sqrt{n}]+1}^{\infty}
			q^j
	=
		2\frac{q^{[\sqrt{n}]+1}}{1-q}
	<
		2\frac{q^{\sqrt{n}}}{1-q}.
\end{equation}

Взявши до уваги оцінки \eqref{r1}, \eqref{r2_2}, \eqref{r3}, \eqref{min|R|} та \eqref{r4}, при $n\geqslant9$ одержимо, що при виконанні умови \eqref{umova_z}
\begin{equation*}
\sum_{k=1}^5
	|\gamma_k(y_0)|
	<
		\frac{212800}{26973n}
		q^{\sqrt{n}}
		\;\frac{1}{1-q^2}
		+
		\frac{140q^{n+\sqrt{n}}}{407n^2}
		+
		\frac
			{40q^{\sqrt{n}}}
			{9(1-q)}
		+
\end{equation*}
\begin{equation*}
		+
		\frac{q}{(1-q)^2}
		\min
			\left\{
				\frac{160}{27(n-\sqrt{n})},
				\frac{8}{3n-7\sqrt{n}}
			\right\}
		+
		\frac
			{2q^{\sqrt{n}}}
			{1-q}
	<
\end{equation*}
\begin{equation*}
	<
		\frac{q^{\sqrt{n}}}{1-q}
		(
			0{,}877
			+
			0{,}0043
			+
			4{,}45
			+
			2
		)
		+
		\frac{q}{(1-q)^2}
		\min
			\left\{
				\frac{160}{27(n-\sqrt{n})},
				\frac{8}{3n-7\sqrt{n}}
			\right\}
	<
\end{equation*}
\begin{equation*}
	<
		\frac{37}{5(1-q)}q^{\sqrt{n}}
		+
		\frac{q}{(1-q)^2}
		\min
			\left\{
				\frac{160}{27(n-\sqrt{n})},
				\frac{8}{3n-7\sqrt{n}}
			\right\}.
\end{equation*}

Лему доведено.

\textbf{5. Доведення теореми 2.}
Відповідно до теореми 5 для встановлення нерівностей \eqref{dno1} і \eqref{dno2} достатньо показати, що для довільних $h>0$, $\beta\in \mathbb{R}$ і всіх номерів $n\geqslant n_h$ ядра $H_{h,\beta}(t)$ задовольняють умову $C_{y_0,2n}$, де $y_0$ --- точка, в якій модуль функції $\Phi_{h,\beta,n}(\cdot)=(H_{h,\beta}\ast\varphi_n)(\cdot)$, $\varphi_n(t)=\textnormal{sign}\sin nt$, досягає найбільшого значення, тобто
\begin{equation*}
|\Phi_{h,\beta,n}(y_0)|=|(H_{h,\beta}\ast\varphi_n)(y_0)|=\|H_{h,\beta}\ast\varphi_n\|_C.
\end{equation*}

Оскільки, як не важко переконатись,
\begin{equation*}
\Phi_{h,\beta,n}(t)
	=
		(H_{h,\beta}\ast\varphi_n)(t)
	=
		\frac{4}{\pi}
		\sum\limits_{\nu=0}^{\infty}
			\frac
				{1}
				{(2\nu+1)\ch((2\nu+1)nh)}
			\sin\left(
				(2\nu+1)nt-\frac{\beta\pi}{2}
			\right),
\end{equation*}
то $\Phi_{h,\beta,n}(\cdot)$ періодична з періодом $2\pi/n$ диференційовна функція і така, що
$\Phi_{h,\beta,n}(\cdot+\frac{\pi}{n})=-\Phi_{h,\beta,n}(\cdot)$. Тому максимальне значення $\pi/n$-періодичної функції $|\Phi_{h,\beta,n}(\cdot)|$ на $[0,\frac{\pi}{n})$ досягається у точці $y_0=y_0(n,h,\beta)=\frac{\theta_n\pi}{n}$, де $\theta_n$ --- єдиний на $[0,1)$ корінь рівняння \eqref{theta}.

Згідно з лемою~2 роботи \cite{My_JAT} для довільного $x\in\mathbb{R}$ і довільного $q\in(0,1)$
\begin{equation}\label{f_x}
\mathcal{P}_q(x)
	>
		\left(
			\frac{1}{2}
			+
			\frac{2q}{(1+q^2)(1-q)}
		\right)
		\left(
			\frac{1-q}{1+q}
		\right)^{\frac {4}{1-q^2}}.
\end{equation}
Тоді з леми 3 і нерівності \eqref{f_x} випливає, що при $n\geqslant9$, $q=e^{-h}$, $k=\overline{1,2n}$, за умов \eqref{umova_n_0} та \eqref{umova_z} виконується нерівність
\begin{equation}\label{geq0}
\mathcal{P}_q(t_k-y_0)
+
\sum\limits_{m=1}^5
	\gamma_m(y_0)
\mathop{\text{sign}} \sin(ny_0-\frac{\beta\pi}{2})
	\geqslant
		0.
\end{equation}
В силу зображення \eqref{SP_Phi_q}, а також нерівностей \eqref{sing_not_zero} і \eqref{geq0} робимо висновок, що при $n\geqslant9$ за умов \eqref{umova_n_0} та \eqref{umova_z} справедливе включення $H_{h,\beta}\in C_{y_0,2n}$. Залишається лише переконатись, що \eqref{umova_z} випливає з \eqref{umova_n_0}.

У роботі \cite{My_zb_2013} було показано, що нерівність \eqref{umova_z} випливає з умови
\begin{equation*}
\frac{43}{10(1-q)}q^{\sqrt{n}}+\frac{q}{(1-q)^2} \min\left\{\frac{160}{57(n-\sqrt{n})}, \frac{8}{3n-7\sqrt{n}}\right\}
 \leqslant
\end{equation*}
\begin{equation}\label{umova_n_00}
\leqslant
\left(\frac{1}{2}+\frac{2q}{(1+q^2)(1-q)}\right)\left(\frac{1-q}{1+q}\right)^{\frac {4}{1-q^2}}.
\end{equation}
При $q=e^{-h}$ безпосередньо переконуємося, що $\eqref{umova_n_0}\Rightarrow\eqref{umova_n_00}$, а отже ${\eqref{umova_n_0}\Rightarrow\eqref{umova_z}}$. 
Теорему доведено.

\textbf{6. Доведення теореми 3.}
Із теорем 1 та 2, а також із співвідношення \eqref{d_m_E_n} випливає, що рівності \eqref{dn} мають місце для усіх номерів $n\geqslant\max\{n_h^*,n_h\}$. Покажемо, що $n_h\geqslant n_h^*$. При $h\geqslant\ln\frac{10}{3}$ вказана нерівність очевидна, оскільки в цьому випадку $n_h^*=1$. Тому залишається переконатись, що при $h\in(0,\ln\frac{10}{3})$ $n_h\geqslant n_h^*=n_h^{**}$.

Покладемо, як і раніше, $q=e^{-h}$. Тоді для доведення теореми достатньо показати, що при $q\in(\frac{3}{10},\,1)$ і $n\geqslant9$ з нерівності
 \begin{equation*}
\frac{37}{5(1-q)}q^{\sqrt{n}}
+
\frac{q}{(1-q)^2} 
\min\left\{
	\frac{160}{27(n-\sqrt{n})}, 
	\frac{8}{3n-7\sqrt{n}}
\right\}
	\leqslant
\end{equation*}
\begin{equation}\label{0umova_n_0}
	\leqslant
		\left(
			\frac{1}{2}+\frac{2q}{(1+q^2)(1-q)}
		\right)
		\left(\frac{1-q}{1+q}\right)^{\frac {4}{1-q^2}}
\end{equation}
випливає нерівність
\begin{equation}\label{0umova_n0}
(1-q)^2
	\geqslant
		\frac{5+3q^2}{1-q^2}
		\frac
			{\left(
				\frac{1+q^2}{2}
			\right)^{2n}}
			{\sqrt{1
				-\left(
					\frac{1+q^2}{2}
				\right)^{2n}}}
		+
		(2+q^{2n})q^{2n}.
\end{equation}

Як доведено в \cite[с.~106]{My_JAT} при $n\geqslant9$ і ${q\in(0,1)}$ із умови $\eqref{umova_n_00}$ випливає нерівність
\begin{equation*}
n
	>
		\frac{160q}{57(1-q)^2} 
		\left(
			\frac{1+q}{1-q}
		\right)^{3}.
\end{equation*}
Оскільки 
\begin{equation*}
\frac{160q}{57(1-q)^2} 
\left(
	\frac{1+q}{1-q}
\right)^{3}
	>
		\frac
			{8q(1+q)}
			{3(1-q)^5}.		
\end{equation*}
то, з урахуванням очевидної імплікації $\eqref{0umova_n_0}\Rightarrow\eqref{umova_n_00}$, одержуємо, що з умови \eqref{0umova_n_0} при $n\geqslant9$ і ${q\in(0,1)}$ випливає нерівність
\begin{equation}\label{ner2}
n
	>
		\frac
			{8q(1+q)}
			{3(1-q)^5}.
\end{equation}

Покажемо, що при ${q\in(\frac{3}{10},1)}$ із \eqref{ner2} випливає нерівність
\begin{equation}\label{ner1}
n
	>
		\frac{5}{1-q^2}
		\ln\frac{2}{1-q}.
\end{equation}

Розглянемо різницю
\begin{equation}\label{vdifference}
v(q)
	=
		\frac
			{8q(1+q)}
			{3(1-q)^5}
		-
		\frac
			{5}
			{1-q^2}
			\ln\frac{2}{1-q},
		\;q\in[\frac{3}{10},1).
\end{equation}
Похідна цієї функції зростає і має вигляд
\begin{equation}\label{vpriveq}
v^\prime(q)
	=
		\frac
			{8(1+6q+3q^2)}
			{3(1-q)^6}
		-
		\frac{5}{(1+q)(1-q)^2}
		\left(
			\frac{2q}{1+q}\ln\frac{2}{1-q}+1
		\right).
\end{equation}
Із відомого розкладу (див., наприклад, \cite[с.~58]{Gradshteyn_1963})
\begin{equation}\label{lnx}
\ln t
	=
		\sum_{k=1}^\infty
			\frac{1}{k}
			\left(
				\frac{t-1}{t}
			\right)^k
		,\;t\geqslant\frac{1}{2},
\end{equation}
при $t=\frac{2}{1-q}$ можемо записати оцінку
\begin{equation}\label{ln21q}
\ln\frac{2}{1-q}
	=
		\sum_{k=1}^\infty
			\frac{1}{k}
			\left(
				\frac{1+q}{2}
			\right)^k
	<
		\sum_{k=1}^\infty
			\left(
				\frac{1+q}{2}
			\right)^k
	=
		\frac{1+q}{1-q}.
\end{equation}
В силу \eqref{vpriveq} і \eqref{ln21q} одержуємо
\begin{equation*}
v^\prime(q)
	>
		\frac
			{8(1+6q+3q^2)}
			{3(1-q)^6}
		-
		\frac{5}{(1-q)^3}.
\end{equation*}

Оскільки $v^\prime(q)$ зростає і $v^\prime(\frac{3}{10})>0$, то $v^\prime(q)>0$, $q\in[\frac{3}{10},\,1)$. Отже, $v(q)$ також зростає на проміжку $[\frac{3}{10},\,1)$, а тому
\begin{equation}\label{vvvvv}
v(q)>v(\frac{3}{10})>0,\; q\in[\frac{3}{10},\,1).
\end{equation}

Із \eqref{vdifference} і \eqref{vvvvv} випливає, що при ${q\in(\frac{3}{10},\;1)}$ $\eqref{ner2}\Rightarrow\eqref{ner1}.$

Нарешті нам залишається довести імплікацію $\eqref{ner1}\Rightarrow\eqref{0umova_n0}.$

Записавши ланцюжок очевидних співвідношень
\begin{equation*}
5\ln
	\frac{2}{1-q}
	>
		\ln
		\frac{16}{(1-q)^5}
	>
		\ln\frac			
			{
				\frac
					{2(5+3q^2)}
					{(1-q^2)(1-q)}
				+3
			}
			{(1-q)^2}
	>
		\ln\frac			
			{
				\frac{5+3q^2}{1-q^2}
				\frac{2}{\sqrt{3(1-q^2)}}
				+3
			}
			{(1-q)^2}
\end{equation*}
і врахувавши нерівність
\begin{equation*}
\ln\frac{2}{1+q^2}
	>
		\frac{1-q^2}{2},
\end{equation*}
яка безпосередньо випливає з розкладу \eqref{lnx} при $t=\frac{2}{1+q^2}$, з \eqref{ner1} отримуємо
\begin{equation*}
n
	>
		\frac
			{1}
			{2\ln
				\frac{1+q^2}{2}
			}
		\ln\frac
			{
				\frac{5+3q^2}{1-q^2}
				\frac{2}{\sqrt{3(1-q^2)}}
				+3
			}
			{(1-q)^2}.
\end{equation*}
Остання нерівність рівносильна наступній нерівності
\begin{equation}\label{umova_n0_0}
(1-q)^2
	>
		\frac{5+3q^2}{1-q^2}
		\frac
			{\left(
				\frac{1+q^2}{2}
			\right)^{2n}}
			{\sqrt{
				1-\left(
					\frac{1+q^2}{2}
				\right)^{2}
				}
			}
		+
		3
		\left(
			\frac{1+q^2}{2}
		\right)^{2n}
\end{equation}
Оскільки $3\left(\frac{1+q^2}{2}\right)^{2n}>(2+q^{2n})q^{2n}$, то із \eqref{umova_n0_0} випливає \eqref{0umova_n0}.
Таким чином
\begin{equation*}
\eqref{0umova_n_0}
	\Rightarrow
		\eqref{ner2}
	\Rightarrow
		\eqref{ner1} 
	\Rightarrow
		\eqref{0umova_n0}.
\end{equation*}

 Теорему доведено.

\textbf{7. Доведення теореми 4.} Знайдемо двосторонні оцінки правої частини формули \eqref{dn}. 
Оскільки
\begin{equation*}
\left|
	\sum\limits_{\nu=1}^{\infty}
		\frac
				{2q^{(2\nu+1)n}}
				{(2\nu+1)(1+q^{2(2\nu+1)n})}
		\sin
		\left(
			(2\nu+1)\theta_n\pi
			-\frac{\beta\pi}{2}
		\right)
\right|
	\leqslant
\end{equation*}
\begin{equation*}
	\leqslant
		\sum\limits_{\nu=1}^{\infty}
			\frac{2q^{(2\nu+1) n}}{(2\nu+1)(1+q^{2(2\nu+1)n})}
	\leqslant
		\frac{2}{3(1+q^{2n})}\frac{q^{3n}}{1-q^{2n}},
\; n\in\mathbb{N},
\end{equation*}
то, враховуючи \eqref{a_n}, одержуємо для довільних $n\in\mathbb{N}$, $q\in(0,1)$ і $\beta\in\mathbb{R}$
\begin{equation*}
\left|
	\sum\limits_{\nu=0}^{\infty}
		\frac
				{2q^{(2\nu+1)n}}
				{(2\nu+1)(1+q^{2(2\nu+1)n})}
		\sin
		\left(
			(2\nu+1)\theta_n\pi
			-\frac{\beta\pi}{2}
		\right)
\right|
	\geqslant
\end{equation*}
\begin{equation*}
	\geqslant
		\frac{2q^n}{1+q^{2n}}
		-\frac{2q^n}{1+q^{2n}}
		\left(
			1
			-|\sin(\theta_n\pi-\frac{\beta\pi}{2})|
		\right)
	-
\end{equation*}
\begin{equation*}
	-
	\left|
		\sum\limits_{\nu=1}^{\infty}
			\frac
				{2q^{(2\nu+1)n}}
				{(2\nu+1)(1+q^{2(2\nu+1)n})}
			\sin
			\left(
				(2\nu+1)\theta_n\pi
				-\frac{\beta\pi}{2}
			\right)
	\right|
\geqslant
\end{equation*}
\begin{equation}\label{ots1}
\geqslant
	\frac{2q^n}{1+q^{2n}}
	\left(
		1
		-\frac{7}{3}
		\frac{q^{2n}}{1-q^{2n}}
	\right),
\end{equation}

\begin{equation*}
\left|
	\sum\limits_{\nu=0}^{\infty}
		\frac
				{2q^{(2\nu+1)n}}
				{(2\nu+1)(1+q^{2(2\nu+1)n})}
		\sin
		\left(
			(2\nu+1)\theta_n\pi
			-\frac{\beta\pi}{2}
		\right)
\right|
	\leqslant
\end{equation*}
\begin{equation*}
	\leqslant
		\frac{2q^n}{1+q^{2n}}
		+\frac{2q^n}{1+q^{2n}}
		\left(
			1
			-|\sin(\theta_n\pi-\frac{\beta\pi}{2})|
		\right)
	+
\end{equation*}
\begin{equation*}
	+
	\left|
		\sum\limits_{\nu=1}^{\infty}
			\frac
				{2q^{(2\nu+1)n}}
				{(2\nu+1)(1+q^{2(2\nu+1)n})}
			\sin
			\left(
				(2\nu+1)\theta_n\pi
				-\frac{\beta\pi}{2}
			\right)
	\right|
\leqslant
\end{equation*}
\begin{equation}\label{ots2}
\leqslant
	\frac{2q^n}{1+q^{2n}}
	\left(
		1
		+\frac{7}{3}\frac{q^{2n}}{1-q^{2n}}
	\right).
\end{equation}
З теореми~3 та оцінок \eqref{ots1} і \eqref{ots2} випливає, що при $n\geqslant n_{h}$ виконується \eqref{Th_3}. 
 Теорему доведено.

\renewcommand{\refname}{}
\makeatletter\renewcommand{\@biblabel}[1]{#1.}\makeatother

\end{document}